%% file: JOGO_indexBB.tex
\documentclass{kluwer} 

\usepackage{latexsym,amssymb}
\usepackage{graphics,times,mathptm,epic,eepic}
\usepackage{epsfig}
\def\beq{\begin{equation}}
\def\eeq{\end{equation}}
\def\bea{\begin{eqnarray}}
\def\eea{\end{eqnarray}}
\def\beas{\begin{eqnarray*}}
\def\eeas{\end{eqnarray*}}

\newcommand{\ds}{\displaystyle}

\newtheorem{lemma}{Lemma}
\newtheorem{theorem}{Theorem}

\begin{document}
\begin{article}
\begin{opening}
\title{Index Branch-and-Bound Algorithm for \\
Lipschitz Univariate Global Optimization with \\
 Multiextremal Constraints\footnotemark[1]}


\author{Yaroslav D. Sergeyev}
\runningtitle{Index Branch-and-Bound Algorithm for  Global
Optimization with Multiextremal Constraints} \institute{ISI-CNR,
c/o DEIS, Universit\`{a} degli Studi della Calabria, Via Pietro
Bucci 41C-42C, 87036 Rende (CS), ITALY and Software Department,
University of Nizhni Novgorod, Gagarin Av. 23, Nizhni Novgorod,
RUSSIAN FEDERATION\\{\sf e-mail} - {\tt yaro@si.deis.unical.it}}

\author{Domenico Famularo}
\institute{DEIS, Universit\`{a} degli Studi della Calabria, Via
Pietro Bucci 41C-42C, 87036 Rende (CS), ITALY\\{\sf e-mail} - {\tt
famularo@deis.unical.it}}

\author{Paolo Pugliese}
\runningauthor{Ya. D. Sergeyev, D. Famularo,  and P. Pugliese}
\institute{DEIS, Universit\`{a} degli Studi della Calabria, Via
Pietro Bucci 41C-42C, 87036 Rende (CS), ITALY\\{\sf e-mail} - {\tt
pugliese@deis.unical.it}}

\footnotetext[1]{ {\bf Acknowledgement.}  The authors thank the
anonymous referees for their great attention to this paper and
very useful and subtle remarks.}

\begin{abstract}
In this paper, Lipschitz univariate constrained global
optimization problems where both the objective function and
constraints  can be multiextremal are considered. The constrained
problem is reduced to a discontinuous unconstrained problem by the
index scheme without introducing additional parameters or
variables. A Branch-and-Bound method that does not use derivatives
for solving the reduced problem is proposed. The method either
determines the infeasibility of the original problem or finds
lower and upper bounds for the global solution. Not all the
constraints are  evaluated during every iteration of the
algorithm, providing a significant acceleration of the search.
Convergence conditions of the new method are established. Test
problems
 and extensive numerical experiments are presented.
 \end{abstract}

\end{opening}

\section{Introduction}
\label{intro} Global optimization problems arise in many real-life
applications and  were intensively studied during last decades
(see, for example,
\cite{Archetti_and_Schoen_(1984),Bomze97,Breiman_and_Cutler_(1993),
Evtushenko_(1992),Floudas_and_Pardalos_(1996),Horst_and_Pardalos_(1995),
Horst_and_Tuy_(1996),Locatelli_and_Schoen_(1999),Lucidi_(1994),
Mladineo_(1992),Pardalos_and_Rosen_(1990),Pinter_(1996),Strongin_(1978),
Sun_and_Li_(1999),Torn_and_Zilinskas_(1989),Zhigljavsky_(1991)},
etc.). Particularly, univariate problems attract attention of many
authors (see
\cite{Calvin_and_Zilinskas_(1999),Hansen_(3),Lamar_(1999),
Locatelli_and_Schoen_(1995),MacLagan_and_Sturge_and_Baritompa_(1996),
Pijavskii_(1972),Sergeyev_(1998),Strongin_(1978),Wang_and_Chang_(1996)})
at least for two reasons. First, there exist a large number of
applications where it is necessary to solve such problems (see
\cite{Brooks_(1958),Hansen_(3),Patwardhan_(1987),Ralston_(1985),Sergeyev_et(1999),Strongin_(1978)}).
Second, there exist numerous schemes (see, for example,
\cite{Floudas_and_Pardalos_(1996),Horst_and_Pardalos_(1995),Horst_and_Tuy_(1996),Mladineo_(1992),Pardalos_and_Rosen_(1990),Pinter_(1996),Strongin_(1978)})
enabling to generalize to the multidimensional case the
mathematical approaches developed to solve univariate problems.

In this paper we consider the global optimization problem
\begin{equation}
\min \{ f(x): x\in [a,b],\hspace{3mm} g_{j}(x)\le 0,\hspace{3mm}  1\le j\le m \},    \label{problem}
\end{equation}
where $f(x)$ and $g_{j}(x), 1\le j\le m,$ are multiextremal
Lipschitz functions (to unify the description process we shall use
the designation $g_{m+1}(x) \triangleq f(x)$). Hereinafter we use
the terminology "multiextremal constraint" to highlight the fact
that the constraints are described by multiextremal functions
$g_{j}(x), 1\le j\le m,$ in the form (\ref{problem}) (of course,
the same subregions of the interval $[a,b]$ may be defined in
another way). In many practical problems the order of the
constraints is fixed and not all the constraints are defined over
the whole search region $[a,b]$ (if the order of the constraints
is not a priori given, the user fixes his/her own ordering in a
way). In the general case, a constraint $g_{j+1}(x)$ is defined
only at subregions where $g_{j}(x)\le 0$. We designate subdomains
of the interval $[a,b]$ corresponding to the set of constraints
from (\ref{problem}) as
\begin{equation}
Q_{1}=[a,b], \hspace{3mm} Q_{j+1}=\{x\in Q_{j}: g_{j}(x)\le 0 \}, \hspace{5mm} 1\le j\le m,  \label{Q}
\end{equation}
\[
Q_{1}\supseteq  Q_{2}\supseteq\ldots\supseteq Q_{m}\supseteq Q_{m+1}.
\]
We introduce the number $M$ such that
\begin{equation}
Q_{M} \neq  \emptyset,  \hspace{5mm} Q_{M+1}=Q_{M+2} \ldots
=Q_{m+1}= \emptyset. \label{M}
\end{equation}
If the feasible region of the problem (\ref{problem}) is not empty
then $Q_{m+1} \not= \emptyset$ and $M=m+1$. In the opposite case
$M$ indicates the last subset $Q_{j}$ from (\ref{Q}) such that
$Q_{j}\neq \emptyset$.

We suppose in this paper that  the functions  $g_{j}(x),1\le j\le
m+1,$ satisfy the Lipschitz condition in the form
\begin{equation}
\mid g_{j}(x') - g_{j}(x'')\mid  \le  L_{j}\mid x'- x''\mid ,
\hspace{3mm} x',x''\in Q_{j},   \hspace{3mm} 1\le j\le m+1.
\label{Lip}
\end{equation}
where the constants
\begin{equation}
0 <  L_{j} < \infty, \hspace{5mm} 1\le j\le m+1 \, , \label{const}
\end{equation}
are known (this  supposition is classical in global optimization
(see \cite{Hansen_(3),Horst_and_Tuy_(1996),Pijavskii_(1972)}), the
problem of estimating the values $L_{j}, 1\le j\le m+1,$ is not
discussed in this paper). Since the functions $g_{j}(x), 1\le
j\le m,$ are supposed to be multiextremal, the subdomains $Q_{j},
2\le j\le m+1,$ can have a few disjoint subregions each. In the
following we shall suppose that all the sets $Q_{j}, 2\le j\le
m+1,$ either are empty or consist of a finite number of disjoint
intervals of a finite positive length.

\begin{figure}[t]
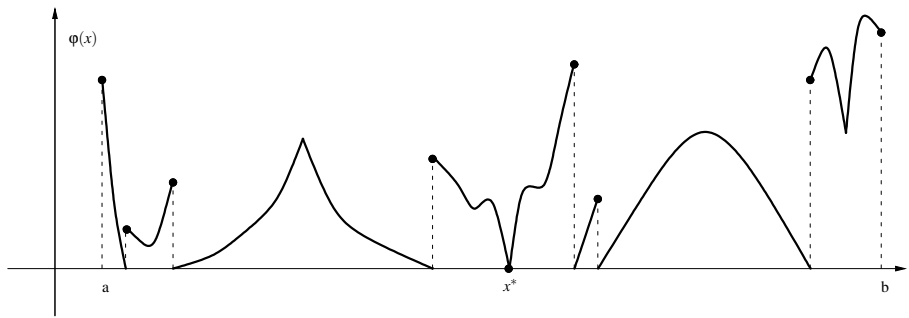

\include{fig1.eepic}
\centerline{a)} \vfill
\include{fig2.eepic}
\centerline{b)} \caption{Construction of the function $\varphi(x)$
\label{figura1}}
\end{figure}

In the example shown in Fig. \ref{figura1},~a) the problem
(\ref{problem}) has two constraints $g_{1}(x)$ and $g_{2}(x)$.
The corresponding sets $Q_{1}=[a,b], Q_{2},$  and $Q_{3}$ are
shown. It can be seen that the subdomain $Q_{2}$ has three
disjoint subregions and the constraint $g_{2}(x)$ is not defined
over the subinterval $[c,d]$. The objective function $f(x)$ is
defined only over the set $Q_{3}$.

The problem (\ref{problem}) may be restated using the
 {\it index scheme} proposed originally in
\cite{Strongin_(1984)} (see also
\cite{Strongin_and_Markin_(1986),Strongin_and_Sergeyev}). The
index scheme does not introduce additional variables and/or
parameters  by opposition to classical approaches in
\cite{Bertsekas_(1996),Bertsekas_(1999),Horst_and_Pardalos_(1995),Horst_and_Tuy_(1996),Nocedal_and_Wright_(1999)}.
It considers constraints one at a time at every point where it has
been decided to calculate $g_{m+1}(x)$. Each constraint $g_i(x)$
is evaluated only if all the inequalities
\[
g_{j}(x)\le 0,  \hspace{5mm} 1\le j<i,
\]
have been satisfied.

In its turn the objective function $g_{m+1}(x)$ is computed only
for that points where all the constraints have been satisfied.

Let us present the index scheme. Using  the designations
(\ref{Q}), (\ref{M}) we can rewrite the problem (\ref{problem})
as the problem of finding a point $x_M^*$ and the corresponding
value $g^{*}_M$ such that
\begin{equation}
g^{*}_M=g_M(x_M^*)=\min\{g_M(x): x\in Q_M\}. \label{6.1.8}
\end{equation}
The values $x_M^*, g^{*}_M$ coincide with the global solution of
the  problem (\ref{problem})  if $M=m+1$, i.e. when the original
problem is feasible. We associate with every point of the interval
$[a,b]$ the {\it index}
\[
\nu = \nu(x),  \hspace{3mm} 1 \le \nu \le M,
\]
which is defined by the conditions
\begin{equation}
g_{j}(x)\le 0,  \hspace{3mm}1\le j\le \nu -1, \hspace{5mm} g_{\nu }(x)>0,       \label{5}
\end{equation}
where for $\nu =m+1$ the last inequality is omitted. We shall
call {\it trial} the operation of evaluation of the functions
$g_{j}(x), 1\le j\le \nu(x),$ at a point $x$. Let us introduce now
an auxiliary function $\varphi (x)$ defined over the interval
$[a,b]$ as follows
\begin{equation}
\varphi (x) = g_{\nu(x) }(x) - \left\{ \begin{array}{ll}
                             0 \, ,          & \mbox{if $\nu (x)<m+1$} \\
                             g^{*}_{m+1} \, ,& \mbox{if $\nu (x)=m+1$}
                                    \end{array}
                            \right.          \label{6}
\end{equation}
where $g^{*}_{m+1}$ is the solution to the problem
(\ref{problem})  and to the problem  (\ref{6.1.8}) in the case
$M=m+1$. Due to (\ref{6.1.8}), (\ref{6}), the function $\varphi
(x)$ has the following properties:
\begin{enumerate}
\item[i.]
$\varphi (x)>0$, when $\nu (x)<m+1$;
\item[ii.]
$\varphi (x)\ge 0$, when $\nu (x)=m+1$;
\item[iii.]
$\varphi (x)=0$, when $\nu (x)=m+1$ and $g_{m+1}(x)=g^{*}_{m+1}$.
\end{enumerate}

In this way the global minimizer of the original constrained
problem (\ref{problem})  coincides with the so\-lu\-tion $x^{*}$
of the following unconstrained discontinuous problem
\begin{equation} \varphi (x^{*})=\min \{ \varphi (x): x\in [a,b]
\},  \label{7} \end{equation} in the case $M=m+1$ and
$g_{m+1}(x^{*})=g^{*}_{m+1}$. Obviously, the value $g^{*}_{m+1}$
used in the construction (\ref{6}) is not known.
Fig.~\ref{figura1}~b) shows the function $\varphi (x)$
constructed for the original problem from Fig.~\ref{figura1}~a).

Numerical methods belonging to the class of {\it information
algorithms} based on probabilistic ideas have been pro\-pos\-ed
for solving  the problem (\ref{7}) in
\cite{Strongin_(1984),Sergeyev_and_Markin_(1995),
Strongin_and_Markin_(1986),Strongin_and_Sergeyev}.

In this paper  a new method called  {\it Index Branch-and-Bound
Algorithm (IBBA)} is introduced for solving  the discontinuous
problem (\ref{7}). The next section shows that, in spite of the
presence of unknown points of discontinuity, it is possible to
construct adaptively improved auxiliary functions (called by the
authors {\it index support functions}) for the function $\varphi
(x)$ and to obtain lower and upper bounds for the global minimum.
The computational scheme of the new method is described in
Section 3. Convergence conditions of the algorithm are
established in Section~4. Section 5 contains wide computational
results showing quite a promising behaviour of the new algorithm.
Finally, Section 6 concludes the paper.

\section{Discontinuous index support functions}
It has been shown in  \cite{Pijavskii_(1972)} that lower and upper
bounds can be found for the global solution $F^*$ of the problem
\begin{equation}
F^* = \min \{ F(x): x\in [a,b]\},    \label{problem1}
\end{equation}
where
\begin{equation}
\mid F(x') - F(x'')\mid  \le  L_F \mid x'- x''\mid , \hspace{5mm}
x',x''\in [a,b],     \label{Lip1}
\end{equation}
through sequential updating of a piece-wise linear support
function
\begin{equation}
\psi (x) \le F(x),  \hspace{5mm} x\in [a,b], \label{support}
\end{equation}
 if the Lipschitz constant $0 <  L_F < \infty \, $ is
known.  The algorithm proposed in \cite{Pijavskii_(1972)}
improves the support function during every iteration by adding a
new point where the objective function $F(x)$ is evaluated. This
procedure enables to draw the support function closer to the
objective and, therefore, to decrease the gap between the lower
and upper bounds. Let us show that by using index approach it is
possible to propose a procedure allowing to obtain lower and
upper bounds for the solution $g^{*}_{m+1}$. In order to induce
the exhaustiveness of the partitioning scheme in the further
consideration it is supposed that constants $K_{j}$ such that
\begin{equation}
L_{j}  < K_{j} < \infty,  \hspace{3mm} 1\le j\le m+1,   \label{K}
\end{equation}
are known. The case $L_{j}  = K_{j}$ is discarded from the
further consideration because in the algorithm of Pijavskii  it
leads to a possibility of generation of a new point coinciding
with one of the points previously generated by the method.

Suppose that $k$ trials have been executed at some points
\begin{equation}
a=x_{0}<x_{1}<\ldots  <x_{i}<\ldots  <x_{k}=b     \label{step1}
\end{equation}
and  the indexes  $\nu_{i}=\nu(x_{i}), 0\le i\le k,$ have been
calculated in accordance with (\ref{5}). Since  the value
$g^{*}_{m+1}$ from (\ref{6}) is not known, it is not possible to
evaluate the function $\varphi (x)$ for the points having the
index $m+1$. In order to overcome this difficulty, we introduce
the function $\varphi_k (x)$ which is evaluated at the points
$x_i$ and gives us the values $z_{i}=\varphi_k (x_i), 0\le i\le
k,$ as follows

\beq
 \varphi_k (x)=g_{\nu (x)}(x) -\left\{
\begin{array}{ll}
                             0          & \mbox{if $\nu (x) < m+1$}\\
                             Z^{*}_{k}      & \mbox{if $\nu(x)=m+1$}

                                    \end{array}
                            \right.            \label{step2.1}
\eeq
 where the value
 \beq
  Z^{*}_{k}=\min \{ g_{m+1}(x_{i}): 0\le i\le k, \nu _{i}=m+1 \}. \label{step2.2}
 \eeq
  estimates $g^{*}_{m+1}$ from (\ref{6}). It can be seen from (\ref{6}), (\ref{step2.1}), and
 (\ref{step2.2}) that  $\varphi_k (x_i)=\varphi
(x_i)$ for all points $x_i$ having indexes $\nu(x_{i})<m+1$ and
\[
0 \le \varphi_k (x_i) \le  \varphi (x_i)
\]
if   $\nu(x_{i})=m+1$.  In addition,
 \beq
  \varphi_k (x) \le 0, \hspace{5mm} x \in \{ x: g_{m+1}(x) \le
  Z^{*}_{k} \}.       \label{step2.3}
 \eeq

During every iteration the trial points $x_{i}, 0\le i\le k,$ form
subintervals
\[
[x_{i-1},x_{i}] \subset [a,b], \hspace{5mm}1\le i\le k,
\]
and every point $x_{i}$ has its own index $\nu_{i}=\nu(x_{i}),
0\le i\le k,$ calculated in accordance with (\ref{5}). Then,
there exist the following three types of subintervals:
\begin{enumerate}
\item[i.] intervals $[x_{i-1},x_{i}]$ such that $\nu_{i-1}=\nu_{i}$;
\item[ii.] intervals $[x_{i-1},x_{i}]$ such that $\nu_{i-1}<\nu_{i}$;
\item[iii.] intervals $[x_{i-1},x_{i}]$ such that $\nu_{i-1}>\nu_{i}$.
\end{enumerate}

The bounding procedure presented below  constructs over each
interval $[x_{i-1},x_{i}]$ for the function $\varphi_k (x)$ from
(\ref{step2.1}) a discontinuos {\it index support function}
$\psi_i (x)$ with the following properties
\[
\psi_i (x) \le \varphi_k (x),\hspace{5mm} x \in [x_{i-1},x_{i}]
\cap Q_{\overline{\nu_{i}}},
\]
where
\[
\overline{\nu_{i}}=\max \{\nu(x_{i-1}),\nu(x_{i}) \}.
\]
Note that the introduced notion is weaker than the usual
definition of a support function (cf. (\ref{support})). In fact,
nothing is required with regard to behaviour of $\psi_i (x)$ over
$[x_{i-1},x_{i}] \setminus Q_{\overline{\nu_{i}}}$ and $\psi_i
(x)$  can be  greater than  $\varphi_k (x)$ on this subdomain.

\begin{figure}[t]
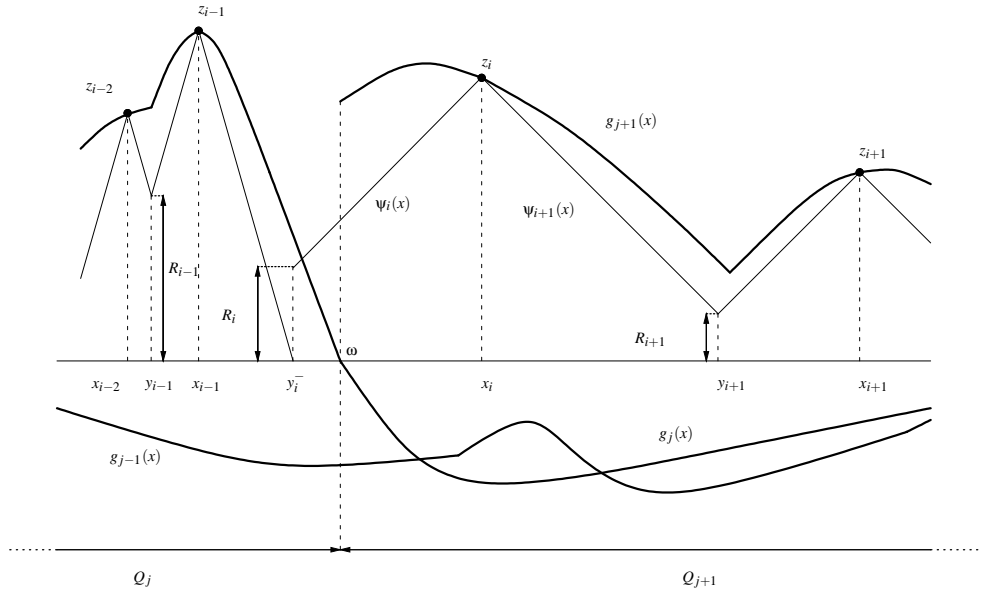

\include{fig3.eepic}
\caption{ The case $\nu(x_{i-2}) = \nu(x_{i-1}) = j$,\hspace{2mm}
$\nu(x_i)=\nu(x_{i+1})=j+1$} \label{figura2}
\end{figure}

Let us consider  one after another  the possibilities (i)-(iii).
The first case, $\nu_{i-1}=\nu_{i}$, is the simplest one. Since
the indexes of the points $x_{i-1},x_{i}$ coincide, the index
support function is similar to that one proposed in
\cite{Pijavskii_(1972)}. In this case, due to  (\ref{Lip}),
(\ref{K}), and \cite{Pijavskii_(1972)}, we can construct for
$\varphi_k (x)$ the index support function $\psi_i(x),x \in
[x_{i-1},x_{i}]$, such  that
\[
\varphi_k (x) \ge \psi_i (x),\hspace{5mm} x \in [x_{i-1},x_{i}]
\cap Q_{\nu_{i}},
\]
where the function $\psi_i (x)$ (see (\ref{step2.1}),
(\ref{step2.2})) has the form
 \beq
  \psi_i (x)=\max \{
g_{\nu_{i}}(x_{i-1}) - K_{\nu_{i}} \mid x_{i-1}-x \mid,
    g_{\nu_{i}}(x_{i}) - K_{\nu_{i}} \mid x_{i}-x \mid \} \label{Pijav_1}
\eeq
 in the case $\nu_{i-1}=\nu_{i}<m+1$ and the form
 \[
 \psi_i (x)=\max \{ g_{m+1}(x_{i-1})-Z^{*}_{k} - K_{m+1} \mid
x_{i-1}-x \mid,
 \]
\beq
    g_{m+1}(x_{i})-Z^{*}_{k} - K_{m+1} \mid x_{i}-x \mid \} \label{Pijav_2}
\eeq
 in the case $\nu_{i-1}=\nu_{i}=m+1$;
 the constants $K_{\nu_{i}}$ are from (\ref {K}). In both cases the global
minimum $R_i$ of the function $\psi_i (x)$ over the interval
$[x_{i-1},x_{i}]$ is
 \beq
R_i=0.5(z_{i-1}+z_{i}-K_{\nu_{i}}(x_{i}-x_{i-1})),  \label {R1}
\eeq
 and is reached at the point
 \beq
y_i=0.5(x_{i-1}+x_{i}-(z_{i}-z_{i-1})/K_{\nu_{i}}).     \label {y}
\eeq
 This case is illustrated in Fig.~\ref{figura2} where the
points $x_{i},x_{i+1}$, ends of the interval $[x_{i},x_{i+1}]$,
have the indexes $\nu_{i}=\nu_{i+1}=j+1<m+1$. In this example
\[
g_{1}(x_{i}) \le 0, \ldots , g_{j}(x_{i}) \le 0, \hspace{5mm} g_{j+1}(x_{i}) > 0,
\]
\[
g_{1}(x_{i+1}) \le 0, \ldots ,g_{j}(x_{i+1}) \le 0, \hspace{5mm}
g_{j+1}(x_{i+1}) > 0,
\]
\[
z_i=\varphi_k (x_i)=g_{j+1}(x_{i}), \hspace{5mm} z_{i+1}=\varphi_k
(x_{i+1})=g_{j+1}(x_{i+1}).
\]
The values $R_{i+1}$ and $y_{i+1}$ are also shown. The interval
$[x_{i-2},x_{i-1}]$ in the same Figure illustrates the case
$\nu_{i-2}= \nu_{i-1}=j$.

The second case is $\nu_{i-1}<\nu_{i}$. Due to the index scheme,
this means that the function $\varphi_k(x)$ has at least one point
of discontinuity $\omega$ over the interval $[x_{i-1},x_{i}]$ (see
an example in Fig.~\ref{figura2}) and consists of parts having
different indexes. To solve the problem (\ref{7}) we are
interested in finding the subregion having the maximal index $M$
from (\ref{M}). The point $x_{i}$ has the index
$\nu_{i}>\nu_{i-1}$ and, due to (\ref{step2.1}), we need an
estimate of the minimal value of the function $\varphi_k(x)$ only
over the domain $[x_{i-1},x_{i}]\cap Q_{\nu_{i}}$. The right
margin of this domain is the point $x_{i}$ because it is the
right end of the interval $[x_{i-1},x_{i}]$ and its index is
equal to $\nu_{i}$. It could be possible to take the point
$x_{i-1}$ as an estimate of the left margin of the domain
$[x_{i-1},x_{i}]\cap Q_{\nu_{i}}$ but a more accurate estimate
can be obtained.

It follows from the inequality $\nu_{i-1}<\nu_{i}$ that
\[
z_{i-1}=\varphi_k (x_{i-1})=g_{\nu_{i-1}}(x_{i-1})>0,\hspace{7mm}
g_{\nu_{i-1}}(x_{i})\le 0.
\]
The function $g_{\nu_{i-1}}(x)$ satisfies the Lipschitz
condition, thus
 \[
 g_{\nu_{i-1}}(x)>0, \hspace{5mm} x \in [x_{i-1},y^{-}_{i})\cap Q_{\nu_{i-1}},
 \]
where the point $y^{-}_{i}$ is obtained from (\ref{Pijav_1})
 \beq
y^{-}_{i}= x_{i-1}+z_{i-1}/K_{\nu_{i-1}}.    \label {y-}
 \eeq
An illustration of this situation is given in Fig.~\ref{figura2}
where the point $\omega \in [y^{-}_{i},x_{i}]$ is such that
$g_{\nu_{i-1}}(\omega)=0 $ and
\[
[x_{i-1},x_{i}]\cap Q_{\nu_{i-1}}=[x_{i-1},\omega],\hspace{5mm}
[x_{i-1},x_{i}]\cap Q_{\nu_{i}}=[\omega,x_{i}],
\]
 \[
[x_{i-1},y^{-}_{i}]\cap Q_{\nu_{i-1}}=[x_{i-1},y^{-}_{i}].
 \]

Therefore, the function $g_{\nu_{i}}(x)$ can be defined at most
over the interval $[y^{-}_{i},x_{i}]$ and the point $y^{-}_{i}$
can be used as an estimate of the left margin of the set
$[x_{i-1},x_{i}]\cap Q_{\nu_{i}}$ for finding a lower bound for
the function $\varphi_k(x)$ over this domain. The corresponding
index support function $\psi_i (x)$ in this case has the form
 \beq
  \psi_i (x)= z_{i} -
K_{\nu_{i}} \mid x_{i}-x \mid  \label {<}
 \eeq
and, therefore,
\[
\min \{ \psi_i (x): x\in [y^{-}_{i},x_{i}] \}\le \min \{ \psi_i
(x): x \in [y^{-}_{i},x_{i}]\cap Q_{\nu_{i}}\}.
\]
This minimum is located at the point $y^{-}_{i}$ and can be
evaluated as
 \beq
 R_i= z_{i}-K_{\nu_{i}}(x_{i}-y^{-}_{i})=
 z_{i}-K_{\nu_{i}}(x_{i}-x_{i-1}-z_{i-1}/K_{\nu_{i-1}}).
\label {R2} \eeq

Let us consider the last case $\nu_{i-1}>\nu_{i}$ being similar to
the previous one. The point $x_{i}$ has the index
$\nu_{i}<\nu_{i-1}$ and, due to the index scheme, we need an
estimate of the minimal value of the function $\varphi_k(x)$ over
the domain $[x_{i-1},x_{i}]\cap Q_{\nu_{i-1}}$.

\begin{figure}[b]
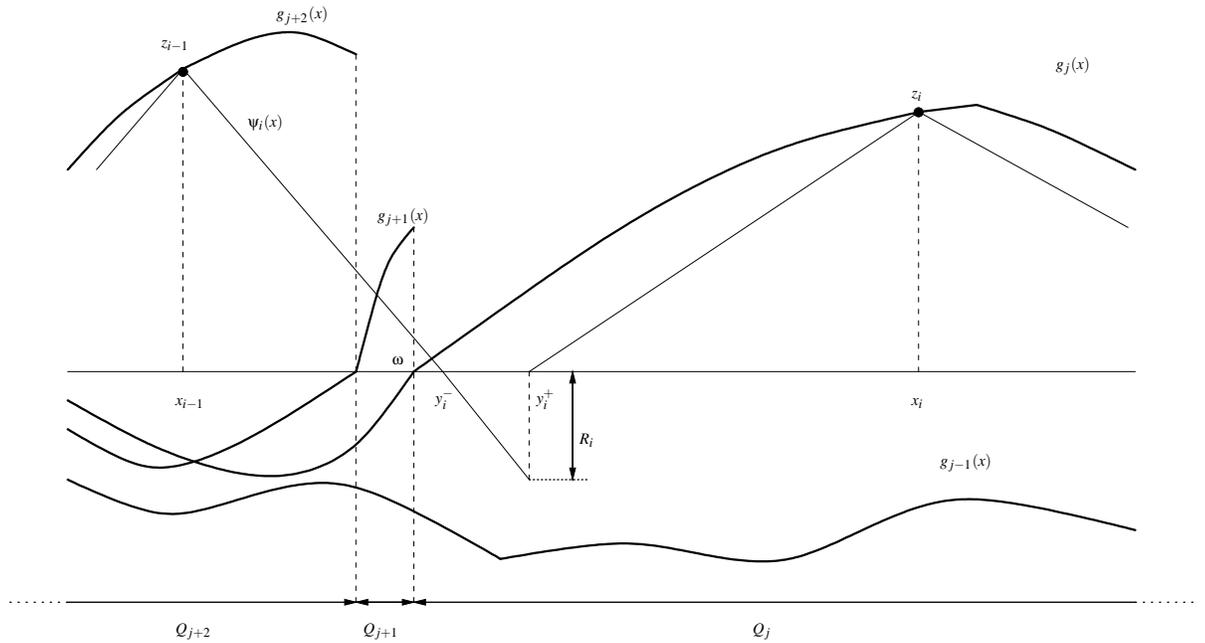

\include{fig4.eepic}
\caption{ The case $\nu(x_{i-1}) = j+2$,\hspace{2mm} $\nu(x_i)=j$}
\label{figura3}
\end{figure}

Since we have $\nu_{i-1}>\nu_{i}$, it follows
\[
z_{i}=\varphi_k (x_{i})=g_{\nu_{i}}(x_{i})>0,\hspace{8mm}
g_{\nu_{i}}(x_{i-1})\le 0.
\]
The function $g_{\nu_{i}}(x)$ satisfies the Lipschitz condition
and, therefore,
\[
g_{\nu_{i}}(x)>0,\hspace{5mm} x\in (y^{+}_{i},x_{i},]\cap
Q_{\nu_{i}},
\]
where
 \beq
  y^{+}_{i}= x_{i}-z_{i}/K_{\nu_{i}}.          \label{y+}
\eeq
 Thus, the function $g_{\nu_{i-1}}(x)$ can be defined at most
over  the interval $[x_{i-1}, y^{+}_{i}]$. The corresponding
index support function
 \beq
  \psi_i (x)= z_{i-1} -
K_{\nu_{i-1}} \mid x_{i-1}-x \mid. \label {>}
 \eeq
 It is evident that
\[
\min \{ \psi_i (x): x\in [x_{i-1}, y^{+}_{i}] \}\le \min \{
\psi_i (x): x\in [x_{i-1}, y^{+}_{i}]\cap Q_{\nu_{i-1}} \}.
\]
It is reached at the point $y^{+}_{i}$ and can be calculated as
 \beq
  R_i= z_{i-1}-K_{\nu_{i-1}}(y^{+}_{i}-x_{i-1})=
z_{i-1}-K_{\nu_{i-1}}(x_{i}-x_{i-1}-z_{i}/K_{\nu_{i}}). \label
{R3}
 \eeq
  This case is illustrated in Fig.~\ref{figura3}. The
points $x_{i-1},x_{i}$ have the indexes $\nu_{i-1}=j+2<m+1,\
\nu_{i}=j$. This means that
\[
g_{1}(x_{i-1}) \le 0, \ldots , g_{j+1}(x_{i-1}) \le 0, \hspace{5mm} g_{j+2}(x_{i-1}) > 0,
\]
\[
g_{1}(x_{i}) \le 0, \ldots , g_{j-1}(x_{i})\le 0, \hspace{5mm} g_{j}(x_{i}) > 0.
\]
The values $z_{i-1}$ and $z_{i}$ are evaluated as follows
\[
z_{i-1}=\varphi_k (x_{i-1})=g_{j+2}(x_{i-1}), \hspace{5mm}
z_{i}=\varphi_k (x_{i})=g_{j}(x_{i}).
\]

Fig.~\ref{figura3} presents a more complex situation in
comparison with Fig.~\ref{figura2}. In fact,
 \[
[x_{i-1},y^{+}_{i}]\cap Q_{\nu_{i-1}}=[x_{i-1},\omega]\setminus
\{Q_{j+1}\cap [x_{i-1},x_{i}]\}.
 \]
The existence of the subregion $Q_{j+1}\cap
[x_{i-1},x_{i}]\neq\varnothing$ cannot be discovered by the
introduced procedure  in the current situation because only the
information
\[
x_{i-1},\nu_{i-1},K_{\nu_{i-1}},z_{i-1},\hspace{5mm} x_{i},\nu_{i},K_{\nu_{i}},z_{i}
\]
regarding the function $\varphi_k (x)$ over $[x_{i-1},x_{i}]$ is
available. This fact is not relevant because we are looking for
subregions with the maximal index $M$, i.e. subregions where the
index is equal to $j+1$ are not of interest because $M \ge j+2$
since $\nu_{i-1}=j+2$.

Now we have completed construction of the function $\psi_i (x)$.
In all three cases, (i) -- (iii), the value $R_i$ being the global
minimum of  $\psi_i (x)$ over the interval $[x_{i-1},x_{i}]$  has
been found (hereinafter we call the value $R_i$ {\it
characteristic} of the interval $[x_{i-1},x_{i}]$). It is
calculated by using one of the formulae (\ref{R1}),(\ref{R2}), or
(\ref{R3}) and is reached at the points $y_{i}$ from (\ref{y}),
$y^{-}_{i}$ is from (\ref{y-}), or $y^{+}_{i}$  from (\ref{y+}),
correspondingly.

If for an interval $[x_{i-1},x_{i}]$ a value $R_i>0$ has been
obtained then, due to the index scheme, it can be concluded that
the global solution $x_{m+1}^*\notin[x_{i-1},x_{i}]$. For
example, in Fig.~\ref{figura2} the intervals
\[
[x_{i-2},x_{i-1}],\hspace{3mm}[x_{i-1},x_{i}],\hspace{3mm}[x_{i},x_{i+1}]
\]
have positive characteristics and, therefore, do not contain the
global minimizer.

Let us now consider an interval $[x_{i-1},x_{i}]$ of the type
(iii) having a negative characteristic (see, for example,  the
interval $[x_{i-1},x_{i}]$ from Fig.~\ref{figura3}). The value
$R_i<0$ has been evaluated at the point $y^{+}_{i}$ as the
minimum of the function $\psi_i (x)$ from (\ref{>}). Since
$z_{i-1}=\psi_i (x_{i-1})>0$ and $R_i=\psi_i (y^{+}_{i})<0$, a
point $\chi\in[x_{i-1},y^{+}_{i}]$ such that $\psi_i (\chi)=0$
can be found. It follows from (\ref{y-}) that $\chi=y^{-}_{i}$.
Thus, the subinterval $[y^{-}_{i},y^{+}_{i}]$ is the only set over
$[x_{i-1},x_{i}]$ where the function $\varphi_k (x)$ can be less
than zero and where, therefore, the global solution $x_{m+1}^*$
can possibly be located. By analogy, it can be shown that when
$R_i<0$ in the cases (i) and (ii), the interval
$[y^{-}_{i},y^{+}_{i}]$ is again the only subinterval of
$[x_{i-1},x_{i}]$  where the global solution  can possibly  be
located.

The {\it Index Branch-and-Bound Algorithm (IBBA)}  proposed in
the next section   at every $(k+1)$th iteration on the basis of
information obtained during the previous $k$ trials constructs
the function $\varphi_k(x)$ and the index support functions
$\psi_i (x),  1\le i\le k$. Among all the intervals
$[x_{i-1},x_{i}],  1\le i\le k,$ it finds an interval $t$ with
the minimal characteristic $R_t$, and chooses the new trial point
$x^{k+1}$ within this interval as follows
 \beq
x^{k+1}=\left\{ \begin{array}{ll}
 0.5(y^{-}_{t}+y^{+}_{t}), &  \hspace{5mm} \nu_{t-1}=\nu_{t}\\
0.5(y^{-}_{t}+x_{t}), &  \hspace{5mm} \nu_{t-1}<\nu_{t}\\
0.5(x_{t-1}+y^{+}_{t}), &  \hspace{5mm} \nu_{t-1}>\nu_{t}
             \end{array}
       \right.            \label{step6}
 \eeq
Note that for intervals having $\nu_{t-1}=\nu_{t}$ the new trial
point $x^{k+1}$ coincides with the Pijavskii point $y_{i}, i=t,$
from (\ref{y}). Thus, the new algorithm at every iteration
updates the function $\varphi_k (x)$ making it closer to $\varphi
(x)$ trying to improve the estimate $Z^{*}_{k}$ of the global
minimum $g^{*}_{m+1}$.

\section{Description of the algorithm}
Let us describe the decision rules of the IBBA. The algorithm
starts with two initial trials at the points $x^{0}=a$ and
$x^{1}=b$. Suppose now that: a search accuracy $\varepsilon$ has
been chosen; $k$ trials have been already done at some points
$x^{0},\ldots ,x^{k}$; their indexes and the value
 \beq
 M^k= \max \{  \nu (x^{i}): 0\le i\le k \}          \label{Mk}
 \eeq
have been calculated. Here the value $M^k$ estimates the maximal
index $M$ from (\ref{M}).

The choice of the point $x^{k+1}, k\ge 1,$ at the $(k+1)$-th
iteration is determined by the rules presented below.

\begin{enumerate}
\item[\bf Step 1.]
The points $x^{0},\ldots  ,x^{k}$ of the previous $k$ iterations
are renumbered by subscripts in order to form the sequence
(\ref{step1}). Thus, two numerations are used during the work of
the algorithm. The record $x^{k}$ means that this point has been
generated during the $k$-th iteration of the IBBA. The record
$x_{k}$ indicates the place of the point in the row (\ref{step1}).
Of course, the second enumeration is changed during every
iteration.

\item[\bf Step 2.]
Recalculate the estimate $Z^{*}_{k}$ from (\ref{step2.2}) and
associate with the points $x_{i}$ the values $z_{i}=\varphi_k
(x_i), 0\le i\le k,$ where the values $\varphi_k (x_{i})$ are from
(\ref{step2.1}).

\item[\bf Step 3.]
For each interval $[x_{i-1},x_{i}], 1\le i\le k$, calculate the
characteristic of the interval
\begin{equation}
R_i=\left\{ \begin{array}{ll}
              0.5(z_{i-1}+z_{i}-K_{\nu_{i}}(x_{i}-x_{i-1})), &
\hspace{5mm} \nu_{i-1}=\nu_{i}\\
              z_{i}-K_{\nu_{i}}(x_{i}-x_{i-1}-z_{i-1}/K_{\nu_{i-1}}), & \hspace{5mm} \nu_{i-1}<\nu_{i} \\
              z_{i-1}-K_{\nu_{i-1}}(x_{i}-x_{i-1}-z_{i}/K_{\nu_{i}}), & \hspace{5mm} \nu_{i-1}>\nu_{i}
             \end{array}
       \right.                  \label{step3}
\end{equation}

\item[\bf Step 4.]
Find the  interval number $t$ such that
\begin{equation}
t = \min \{ \arg \min \{ R_i: 1\le i\le k \} \}.   \label{step4}
\end{equation}

\item[{\bf Step 5.}] {\it (Stopping Rule)}
If $R_t>0$, then Stop (the feasible region is empty). Otherwise, if
 \begin{equation}
x_{t}-x_{t-1} > \varepsilon             \label{step5}
\end{equation}
go to Step 6 ($\varepsilon$ is a preset accuracy and $t$ is from
(\ref{step4})). In the opposite case, Stop (the required accuracy
has been reached).

\item[\bf Step 6.]
Execute the $(k+1)$-th trial at the point $x^{k+1}$ from
(\ref{step6}),  evaluate its index $\nu (x^{k+1})$ and the
estimate $M^{k+1}$, and go to Step~1.
\end{enumerate}

In the following section we will gain more insight the method by
establishing and discussing its convergence conditions.

\section{Convergence conditions}

In this section we demonstrate that the infinite trial sequence
$\{x^k\}$ generated by the algorithm IBBA ($\varepsilon=0$ in the
stopping rule) converges to the global solution of the
unconstrained problem (\ref{7}) and, as consequence, to the
global solution of the initial constrained problem
(\ref{problem})  if it is feasible. In the opposite case the
method establishes infeasibility of the problem (\ref{problem})
in a finite number of iterations.

In Lemma 1, we prove the exhaustiveness of the branching scheme.
The convergence results of the proposed method can be derived as
a particular case of  general convergence studies given in
\cite{Horst_and_Tuy_(1996),Pinter_(1996),Sergeyev_(1999)}. We
present a detailed and independent proof of these results in
Theorems 1 and 2.

\begin{lemma}
Let $\bar{x}$  be a limit point of the sequence $\{x^{k}\}$
generated by the IBBA with $\varepsilon = 0$ in the stopping rule
(\ref{step5}), and let $i=i(k)$ be the number of an interval
$[x_{i(k)-1},x_{i(k)}]$ containing this point during the $k$-th
iteration. Then
\begin{equation}
\lim_{k \rightarrow \infty} x_{i(k)}-x_{i(k)-1}= 0.     \label{lim}
\end{equation}
\end{lemma}

\noindent {\bf Proof:} During the current $k$-th iteration an
interval $[x_{t-1},x_{t}]$  is chosen for subdivision. Due to the
decision rules of the IBBA and (\ref{y-}), (\ref{y+}), this means
that its characteristic $R_t\le0$ and the point $x^{k+1}$ from
(\ref{step6}) falling into the interval $(x_{t-1},x_{t})$  can be
rewritten as follows
 \beq
 x^{k+1}=\left\{ \begin{array}{ll}
 0.5(x_{t}+x_{t-1}+z_{t-1}/K_{\nu_{t-1}}-z_{t}/K_{\nu_{t}}), &  \hspace{5mm} \nu_{t-1}=\nu_{t}\\
  0.5(x_{t}+x_{t-1}+z_{t-1}/K_{\nu_{t-1}}), &  \hspace{5mm} \nu_{t-1}<\nu_{t}\\
   0.5(x_{t}+x_{t-1}-z_{t}/K_{\nu_{t}}), &  \hspace{5mm} \nu_{t-1}>\nu_{t}\\
  \end{array}
       \right.   \label{x_explcit}
 \eeq
This point divides the interval $[x_{t-1},x_{t}]$ into two
subintervals
\begin{equation}
[x_{t-1},x^{k+1}], \hspace{5mm} [x^{k+1},x_{t}].            \label{interv}
\end{equation}
Let us show that the following contracting estimate
\[
\max \{ x_{t}-x^{k+1}, x^{k+1}-x_{t-1} \} \le
\]
\begin{equation}
0.5 (1+\max \{  L_{\nu_{t-1}}/K_{\nu_{t-1}}, L_{\nu_{t}}/K_{\nu_{t}}\})(x_{t}-x_{t-1})      \label{contract1}
\end{equation}
holds for the intervals (\ref{interv}), where
\begin{equation}
0.5 \le 0.5 (1+\max \{  L_{\nu_{t-1}}/K_{\nu_{t-1}}, L_{\nu_{t}}/K_{\nu_{t}}\}) < 1.   \label{contract2}
\end{equation}
Let us consider three cases.

i. In the first case $\nu_{t-1}=\nu_{t}$, and $z_{t-1}=g_{\nu_{t}}
(x_{t-1})$, $z_t=g_{\nu_{t}} (x_t)$. It follows from (\ref{Lip1}),
(\ref{K}) that
\[
 L_{\nu_{t-1}}=L_{\nu_{t}}< K_{\nu_{t-1}}=K_{\nu_{t}},
\]
\[
\mid z_t-z_{t-1} \mid \leq L_{\nu_{t}} (x_t-x_{t-1})
 \]
 and, due to (\ref{y}), we have
\[
\max \{x_{t}-x^{k+1}, x^{k+1}-x_{t-1} \} \le 0.5 (1+
L_{\nu_{t-1}}/K_{\nu_{t-1}})(x_{t}-x_{t-1}).
\]
Thus, (\ref{contract1}) and (\ref{contract2}) have been
established.

ii. In the second case, $\nu_{t-1}<\nu_{t}$, and, therefore, due
to the index scheme, $z_{t-1}=g_{\nu_{t-1}}(x_{t-1})>0$ and
$g_{\nu_{t-1}}(x_{t})\le 0$. From this estimate and the obvious
relation
\[
g_{\nu_{t-1}}(x_{t})\ge z_{t-1}- L_{\nu_{t-1}}(x_{t}-x_{t-1})
\]
we obtain
 \beq
 z_{t-1}- L_{\nu_{t-1}}(x_{t}-x_{t-1}) \le 0. \label{key}
\eeq
 Since $z_{t-1} > 0$, it follows from (\ref{key}), (\ref{step6}),
  and (\ref{K}) that
\[
x^{k+1}-x_{t-1}= 0.5 (x_{t}-x_{t-1}+ z_{t-1}/K_{\nu_{t-1}}) \le
\]
\beq
 0.5 (x_{t}-x_{t-1}+
L_{\nu_{t-1}}/K_{\nu_{t-1}}(x_{t}-x_{t-1})). \label{contract1.1}
\eeq
 Let us now estimate the difference $x_{t}-x^{k+1}$.
 \beq
  x_{t}-x^{k+1} = 0.5 (x_{t}-x_{t-1}-
z_{t-1}/K_{\nu_{t-1}}) < 0.5 (x_{t}-x_{t-1}). \label{contract1.2}
 \eeq
Obviously, the estimate (\ref{contract1}) is the result of
(\ref{contract1.1}) and (\ref{contract1.2}).

iii. The case $\nu_{t-1}>\nu_{t}$ is considered by a complete
analogy to the case (ii) and leads to  estimates
 \[
x^{k+1}-x_{t-1} < 0.5 (x_{t}-x_{t-1}),
\]
\[
 x_{t}-x^{k+1} \le 0.5 (x_{t}-x_{t-1}+
L_{\nu_{t}}/K_{\nu_{t}}(x_{t}-x_{t-1})).
 \]

To prove (\ref{contract2}) it is enough to mention that
$L_{\nu_{t-1}},K_{\nu_{t-1}},L_{\nu_{t}},$ and $K_{\nu_{t}}$ are
constants and (\ref{K}) takes place for them. The result
(\ref{lim}) is a straightforward consequence of the decision
rules of the IBBA and the estimates (\ref{contract1}),
(\ref{contract2}). \rule{5pt}{5pt}

\begin{theorem} If the original problem (\ref{problem})  is infeasible then the algorithm
stops in a finite number of iterations.
\end{theorem}

\noindent {\bf Proof:} If the original problem (\ref{problem}),
(\ref{Lip}) is infeasible then the maximal index $M$ over the
interval $[a,b]$ is less than $m+1$. In this case (see (\ref{5}),
(\ref{6}), and (\ref{step2.1}))
\[
\varphi (x)= \varphi_k (x)> 0, \hspace{5mm} x \in [a,b].
\]
On one hand, due to (\ref{Lip}), (\ref{K}), the linear pieces of
the index support functions $\psi_i (x), 1\le i\le k,$ from
(\ref{Pijav_1}),  (\ref{<}), and (\ref{>}) constructed by the
algorithm have a finite slope. On the other hand, Lemma 1 shows
that the length of any interval containing any limit point goes to
zero.

Thus, it follows from our supposition regarding the sets $Q_{j},
1\le j\le m+1,$  being either empty or consisting of a finite
number of disjoint  intervals of a finite positive length and the
formulae (\ref{contract2}), (\ref{step3}), and (\ref{step4}) that
there exists a finite iteration number $N$ such that a
characteristic $R_{t(N)}>0$ will be obtained and the algorithm
will stop. \rule{5pt}{5pt}

Let us now consider the case when the original problem
(\ref{problem})  is feasible. This means that $M=m+1$ in
(\ref{6.1.8}), (\ref{6}). Let us denote by $X^*$ the set of the
global minimizers of the problem (\ref{problem})  and by $X'$ the
set of limit points of the sequence $\{x^{k}\}$ generated by the
IBBA with $\varepsilon = 0$ in the stopping rule (\ref{step5}).

\begin{theorem}
If the problem (\ref{problem})  is feasible then $X^*=X'$.
\end{theorem}

\noindent {\bf Proof:} Since the problem (\ref{problem})  is
feasible, the sets $Q_{j}, 1\le j\le m+1,$  are not empty and
therefore, due to our hypotheses, they consist of a finite number
of disjoint intervals of a finite positive length. This fact
together with $\varepsilon$ from (\ref{step5}) equal to zero
leads to  existence of an iteration $q$ during which a point
$x^{q}$  having the index $m+1$ will be generated. Thus, (see
(\ref{step2.1}), (\ref{step2.2})) the first value $z_i=0$
corresponding to the point $x^{q}$ will be obtained and during
all the iterations $k>q$ there will exist at least two intervals
having negative characteristics (see (\ref{step3})).

Let us return to the interval $[x_{i-1},x_{i}]$ from Lemma 1
containing a limit point $\bar{x} \in X'$. Since it contains the
limit point and the trial points are chosen by the rule
(\ref{step4}), its characteristic should be negative too for all
iterations $k>q$. Then, by taking into consideration the facts
that  $z_i \ge 0, 1 \le i \le k,$ (see (\ref{step2.1})) it follows
from Lemma~1 and (\ref{step3}), (\ref{step4}) that
 \beq
 \lim_{k \rightarrow \infty}{R_{i(k)}} = 0.  \label{x'}
\eeq
 We can conclude from (\ref{x'}) and $R_{i(k)}<0,k>q,$ that
\beq
 \lim_{k \rightarrow \infty}{\varphi_k (\bar{x})}=\varphi (\bar{x})=
 0.   \label{varphix'}
\eeq

Let us consider an  interval $[x_{j(k)-1},x_{j(k)}]$ containing a
global minimizer $x^{*} \in X^*$ during an iteration $k>q$. At
first, we show that there will exist an iteration number $c \ge
q$ such that $\nu _{j(c)-1}=m+1$ or $\nu _{j(c)}=m+1$. If trials
will fall within the interval $[x_{j(k)-1},x_{j(k)}]$, due to the
decision rules of the IBBA, such a trial will be generated.
Suppose that trials will not fall into this interval and
\[
\Gamma=\max \{ \nu(x_{j(k)-1}), \nu(x_{j(k)}) \}<m+1.
 \]
  The point
$x^{*}$ is feasible, this means that
\[
x^{*} \in [\alpha,\beta]=[x_{j(k)-1},x_{j(k)}]\cap Q_{m+1},
\]
where the interval $[\alpha,\beta]$ has a finite positive length
and
 \[
g_{l}(x)\le 0,  \hspace{5mm} 1\le l\le m, \hspace{5mm} x \in
[\alpha,\beta].
  \]
We obtain from these inequalities that, due to (\ref{step3}) and
(\ref{K}), the characteristic
 \beq
  R_{j(k)} \le \min \{
g_{\Gamma}(x): x \in [\alpha,\beta] \}<0. \label{gamma}
 \eeq
Since trials do not fall at the interval $[x_{j(k)-1},x_{j(k)}]$,
it follows from (\ref{step3}) that $R_{j(k)}$  is not changed from
iteration to iteration. On the other hand, the characteristic
$R_{i(k)}\rightarrow 0$ when $k \rightarrow \infty$. This means
that at an iteration number $k'>q$ the characteristic of the
interval $[x_{i-1},x_{i}], i = i(k')$, will not be minimal. Thus,
a trial will fall into the interval $[x_{j-1},x_{j}]$. The
obtained contradiction proves generation of a point
 \[
 x^c\in [\alpha,\beta],\hspace{5mm} c\ge k',\hspace{2mm}  \nu(x^c)=m+1.
 \]

We can now estimate the characteristic $R_{j(k)}$ of the interval
$[x_{j(k)-1},x_{j(k)}]$ containing the global minimizer $x^{*} \in
X^*$ during an iteration $k>c$. We have shown that at least one of
the points $x_{j(k)-1},x_{j(k)}$ will have the index $m+1$. Again,
three cases can be considered.

In the case $\nu _{j-1}=\nu _{j}=m+1$ it follows from
(\ref{Lip}), (\ref{step2.1}) that
\[
z_{j-1}-\varphi_k (x^{*})\le  L_{m+1}(x^{*}-x_{j-1}).
\]
From (\ref{step2.3}) we have $\varphi_k (x^{*})\le 0$ and,
therefore,
\[
z_{j-1}\le  L_{m+1}(x^{*}-x_{j-1}),
\]
Analogously, for the value $z_{j}$ it follows
\[
z_{j}\le  L_{m+1}(x_{j}-x^{*}).
\]
From these two estimates we obtain
\[
z_{j}+z_{j-1}\le  L_{m+1}(x_{j}-x_{j-1}).
\]
By using (\ref{K}), (\ref{R1}),  and the last inequality we deduce
\begin{equation}
R_{j(k)} \le  (L_{m+1}-K_{m+1})(x_{j(k)}-x_{j(k)-1}) < 0.   \label{R*1}
\end{equation}
Analogously, it can be seen from (\ref{y-}) and (\ref{R2})  that
in the case $\nu _{j-1}< \nu _{j}$ the estimate
\begin{equation}
R_{j(k)} \le  (L_{m+1}-K_{m+1})(x_{j(k)}-y^{-}_{j(k)}) < 0
\label{R*2}
\end{equation}
takes place because
\[
z_{j}\le  L_{m+1}(x_{j}-x^{*}) \le L_{m+1}(x_{j(k)}-y^{-}_{j(k)}).
\]
For the case $\nu _{j-1}> \nu _{j}$ (see (\ref{y+}) and
(\ref{R3})) we have
\begin{equation}
R_{j(k)} \le  (L_{m+1}-K_{m+1})(y^{+}_{j(k)}-x_{j(k)-1}) < 0.   \label{R*3}
\end{equation}

It follows from (\ref{gamma}) -- (\ref{R*3})
that the characteristic of the interval
$[x_{j-1},x_{j}]$ containing the global minimizer $x^{*}$ will be
always negative. Assume now, that $x^{*}$ is not a limit point of
the sequence $\{x^{k}\}$, then there exists a number $P$ such
that for all $k \ge  P$ the interval $[x_{j-1},x_{j}],j=j(k)$, is
not changed, i.e. new points will not fall into this interval
and, as a consequence, its characteristic $R_{j(k)}$ will not
change too.

Consider again the interval $[x_{i-1},x_{i}]$ from Lemma 1
containing a limit point $\bar{x} \in X'$. It follows from
(\ref{x'}) and the fact that $R_{j(k)}$ is a negative constant
that there exists an iteration number $N$ such that
\[
R_{j(N)} < R_{i(N)}.
\]
Due to decision rules of the IBBA, this means that a trial will
fall into the interval $[x_{j-1},x_{j}]$. But this fact
contradicts  our assumption that $x^{*}$ is not a limit point.

Suppose now that there exists a limit point $\bar{x} \in X'$ such
that $\bar{x} \not\in X^{*}$. This means that $  \varphi (\bar{x})
> \varphi (x^{*}), x^{*} \in X^{*}$. Impossibility of this fact
comes from (\ref{x'}), (\ref{varphix'}), and the fact of $x^{*}
\in X'$. \rule{5pt}{5pt}

We can conclude that if the algorithm has stopped and has not
established that $Q_{m+1}= \emptyset$ then the following
situations are possible:
\begin{enumerate}
\item[i.]
If $M^k<m+1$, then this means that the accuracy $\varepsilon$ was
not sufficient for establishing the feasibility of the  problem;
\item[ii.]
If $M^k=m+1$ and all the intervals $[x_{p-1},x_{p}]$ such that
\beq
 \max \{\nu_{p-1}, \nu_p \} < m+1          \label{maxp}
 \eeq
 have positive  characteristics
then, we can conclude that the global minimum $z^{*}$ of the
original problem (\ref{problem}) can be bounded as follows
\[
z^{*} \in [R_{t(k)}+Z^{*}_{k}, Z^{*}_{k}],
\]
where the value $Z^{*}_{k}$ is from (\ref{step2.2}) and $R_{t(k)}$
is the characteristic corresponding to the interval  number
$t=t(k)$ from (\ref{step4}).
\item[iii.]
If $M^k=m+1$ and there exists an interval $[x_{p-1},x_{p}]$ such
that $R_p \le 0$ and   (\ref{maxp}) takes place then, the value
$Z^{*}_{k}$ can be taken as an upper bound of the  global minimum
$z^{*}$ of the original problem (\ref{problem}). A rouge lower
bound can be calculated easily by taking the trial points $x_i$
such that $\nu(x_i)= m+1$ and constructing for $f(x)$  the
support function of the type \cite{Pijavskii_(1972)} using only
these points. The global minimum of this support function over
the search region $[a,b]$ will be a lower bound for $z^{*}$. A
more precise lower bound can be obtained by minimizing this
support function over the set
\[
\bigcup[x_{i-1},x_{i}], \hspace{5mm}  R_i<0,\hspace{3mm} 1\le
i\le k.
\]
\end{enumerate}

We do not discuss here the peculiarities of the implementation of
the IBBA. Let us make only two remarks. First, it is not necessary
to re-calculate all the characteristics during Step~3 but it is
sufficient to do this operation only for two new intervals
generated during the previous iteration. Second, as it follows
from the proofs of Theorems~$1, 2$, it is possible to exclude from
consideration all the intervals having positive characteristics.

\section{Numerical comparison}
The IBBA algorithm has been numerically compared to the method
(indicated hereinafter as {\it PEN}) proposed by Pijavskii (see
\cite{Pijavskii_(1972),Hansen_(3)}) combined with a  penalty
function. The PEN has been chosen for comparison because it uses
the same information about the problem as the IBBA -- the
Lipschitz constants for the objective function and constraints.

Ten differentiable and ten non-differentiable test problems
introduced in \cite{Famularo_Sergeyev_Pugliese} have been used.
In addition, the IBBA has been applied to one differentiable and
one non-differentiable infeasible test problem from
\cite{Famularo_Sergeyev_Pugliese}.  Since the order of constraints
can influence speed of the IBBA significantly, it has  been chosen
the same as in \cite{Famularo_Sergeyev_Pugliese}, without
determining the best order for the IBBA. The same accuracy $
\varepsilon = 10^{-4} \left(b - a\right)$ (where $b$ and $a$ are
from (\ref{problem})) has been used in all the experiments for
both methods.

\begin{center}
\begin{table}[t]
\caption{Results of the experiments executed by the IBBA with the
differentiable problems.} \label{tabpr:1}
\begin{tabular}{l|rrrrrrrr} \hline\noalign{\smallskip}
Problem & XIBBA & FIBBA & $N_{g_1}$ & $N_{g_2}$ & $N_{g_3}$ &
$N_{f}$ & Iterations & Eval. \\
\noalign{\smallskip}\hline\noalign{\smallskip} 1 & $1.05726259$ &
$-7.61226549$ & $10$ & $-$ & $-$ & $13$ & $23$ & $36$ \\
 2 &
$1.01559921$ & $5.46160556$ & $206$ & $-$ & $-$ & $21$ & $227$ &
$248$  \\
 3 & $-5.99182849$ & $-2.94266082$ & $40$ & $-$ & $-$ &
$22$ & $62$ & $84$ \\
 4 & $2.45965829$ & $2.84080900$ & $622$ &
$156$ & $-$ & $175$ & $953$ & $1459$ \\
5 & $9.28501542$ & $-1.27484676$ & $8$ & $14$ & $-$ & $122$ &
$144$ & $402$ \\
 6 & $2.32396546$ & $-1.68515824$ & $14$ & $80$ &
$-$ & $18$ & $112$ & $228$  \\
 7 & $-0.77473979$ & $-0.33007410$ &
$35$ & $18$ & $-$ & $241$ & $294$ & $794$ \\
 8 & $-1.12721979$ &
$-6.60059664$ & $107$ & $43$ & $5$ & $82$ & $237$ & $536$ \\
 9 &
$4.00046339$ & $1.92220990$ & $7$ & $36$ & $6$ & $51$ & $100$ &
$301$  \\
10 &$4.22474504$ & $1.47400000$ & $37$ & $15$ & $195$ & $1173$ &
$1420$ & $5344$  \\ Average & $-$ & $-$ &
$-$ & $-$ & $-$ & $-$ & $357.2$ & $943.2$ \\
\noalign{\smallskip}\hline
\end{tabular}
\end{table}
\end{center}

\begin{center}
\begin{table}[h]
\caption{Results of the experiments executed by the IBBA with the
non-differentiable problems.\label{nondif}}
    \begin{tabular}{l|rrrrrrrr} \hline\noalign{\smallskip}
Problem & XIBBA & FIBBA & $N_{g_1}$ & $N_{g_2}$ & $N_{g_3}$ &
$N_{f}$ & Iterations & Eval. \\
\noalign{\smallskip}\hline\noalign{\smallskip}
 1 & $1.25830963$ & $4.17420017$ & $23$ & $-$ & $-$ & $28$ &   $51$ & $79$ \\
 2 & $1.95967593$ & $-0.07915191$ & $18$  & $-$ & $-$ & $16$ & $34$ & $50$   \\
 3 & $9.40068508$ & $-4.40068508$ & $171$ & $-$ & $-$ & $19$ & $190$ & $209$   \\
 4 & $0.33286804$ & $3.34619770$ & $136$ & $15$ & $-$ & $84$ & $235$ & $418$   \\
 5 & $0.86995142$ & $0.74167893$ & $168$ & $91$ & $-$ & $24$ & $283$ & $422$  \\
 6 & $3.76977185$ & $0.16666667$ & $16$ & $16$ & $-$ & $597$ & $629$ & $1839$ \\
 7 & $5.20120767$ & $0.90312158$ & $63$ & $18$ & $-$ & $39$ &  $120$ & $216$ \\
 8 & $8.02835033$ & $4.05006890$ & $29$ & $11$ & $3$ & $21$ &  $64$ & $144$ \\
 9 & $0.95032461$ & $2.64804102$ & $8$ & $86$ & $57$ & $183$ & $334$ & $1083$ \\
10 & $0.79996352$ & $1.00023345$ & $42$ & $3$ & $17$ & $13$ &  $75$ & $151$ \\
 Average & $-$ & $-$ & $-$ & $-$ & $-$ & $-$ & $201.5$ &
$461.1$ \\ \noalign{\smallskip}\hline
    \end{tabular}
\end{table}
    \end{center}

\vspace{-20mm}

In Table \ref{tabpr:1} (Differentiable problems) and Table II
(Non-Differentiable problems) the results obtained by the IBBA
have been summarized and the columns in the Tables have the
following meaning:
\begin{enumerate}
\item [-] the columns XIBBA and FIBBA represent the estimate   to the global
 solution $(x^*,f(x^*))$  found by the IBBA for each problem;
\item [-] the columns $N_{g_1}$, $N_{g_2}$,  $N_{g_3}$   represent
 the number of trials where the constraint $g_i, 1 \le i \le 3,$ was the last
 evaluated constraint;
\item [-] the column $N_{f}$ shows how many times the objective function $f(x)$
has been evaluated;
\item [-] the column "Eval." is the total number of evaluations
 of the objective function and the constraints. This quantity is equal
to:
\begin{enumerate}
\item [-]  $N_{g_1} + 2 \times N_{f} ,$\,for problems with one
constraint;
\item [-] $ N_{g_1} + 2 \times N_{g_2} + 3 \times N_{f},$\,  for problems with two
constraints;
\item [-] $ N_{g_1} + 2 \times N_{g_2} + 3 \times N_{g_3} + 4 \times N_{f},$\, for  problems with three
constraints.
\end{enumerate}
\end{enumerate}

\begin{center}
\begin{table}[t]
\caption{Differentiable functions. Numerical results obtained by
the PEN.}
    \label{pentab:1}
    \begin{tabular}{l|rrrrr} \hline\noalign{\smallskip}
Problem & XPEN & FXPEN & $P^*$ & Iterations & Eval. \\
\noalign{\smallskip}\hline\noalign{\smallskip} 1 & $1.05718004$ &
$-7.61185807$ & $15$ & $83$ & $166$ \\ 2 & $1.01609254$ &
$5.46142698$ & $90$ & $954$ & $1906$ \\ 3 & $-5.99184997$ &
$-2.94292577$ & $15$ & $119$ & $238$ \\ 4 & $2.45953057$ &
$2.84080890$ & $490$ & $1762$ & $5286$ \\ 5 & $9.28468704$ &
$-1.27484673$ & $15$ & $765$ & $2295$ \\ 6 & $2.32334492$ &
$-1.68307049$ & $15$ & $477$ & $1431$ \\ 7 & $-0.77476915$ &
$-0.33007412$ & $15$ & $917$ & $2751$ \\ 8 & $-1.12719146$ &
$-6.60059658$ & $15$ & $821$ & $3284$ \\ 9 & $4.00042801$ &
$1.92220821$ & $15$ & $262$ & $1048$ \\ 10 & $4.22482084$ &
$1.47400000$ & $15$ & $2019$ & $8076$ \\ Average & $-$ & $-$ & $-$
& $817.9$ & $2648.1$  \\
\noalign{\smallskip}\hline\noalign{\smallskip}
\end{tabular}
\end{table}
\end{center}

\begin{center}
\begin{table}
\caption{Non-Differentiable problems. Numerical results obtained
by the PEN.}
    \label{pentab:2}
    \begin{tabular}{l|rrrrr} \hline\noalign{\smallskip}
Problem & XPEN & FXPEN & $P^*$ & Iterations & Eval. \\
\noalign{\smallskip}\hline\noalign{\smallskip} 1 & $1.25810384$ &
$4.17441502$ & $15$ & $247$ & $494$ \\ 2 & $1.95953624$ &
$-0.07902265$ & $15$ & $241$ & $482$ \\ 3 & $9.40072023$ &
$-4.40072023$ & $15$ & $797$ & $1594$ \\ 4 & $0.33278550$ &
$3.34620350$ & $15$ & $272$ & $819$ \\ 5 & $0.86995489$ &
$0.74168456$ & $20$ & $671$ & $2013$ \\ 6 & $3.76944805$ &
$0.16666667$ & $15$ & $909$ & $2727$ \\ 7 & $5.20113260$ &
$0.90351752$ & $15$ & $199$ & $597$ \\ 8 & $8.02859874$ &
$4.05157770$ & $15$ & $365$ & $1460$ \\ 9 & $0.95019236$ &
$2.64804101$ & $15$ & $1183$ & $4732$ \\ 10 & $0.79988668$ &
$1.00072517$ & $15$ & $135$ & $540$ \\ Average & $-$ & $-$ & $-$ &
$501.9$ & $1545.8$ \\
\noalign{\smallskip}\hline\noalign{\smallskip}
    \end{tabular}
\end{table}
    \end{center}

In Table \ref{pentab:1} (Differentiable problems) and Table
\ref{pentab:2} (Non-Differentiable problems) the results obtained
by the  PEN are collected. The constrained problems were reduced
to the unconstrained ones as follows

\beq f_{P^*}(x)=f(x) + P^*
\max\left\{g_1(x),g_2(x),\dots,g_{N_v}(x),0\right\}.  \label{pen}
\eeq
The coefficient $P^*$ has been computed by the rules:
\begin{enumerate}
\item the coefficient $P^*$ has been  chosen equal to $15$ for all  the problems and it has  been checked if the found solution (XPEN,FXPEN) for each problem belongs or not to the feasible subregions;
\item if it does not belong to the feasible subregions, the coefficient $P^*$ has been iteratively increased by $10$ starting from $20$ until a feasible solution has been found. Particularly, this means that a feasible solution has not been found in Table \ref{pentab:1} for the problem 2 when $P^*$ is equal to $80$, for the problem 4 when $P^*$ is equal to $480$, and in Table \ref{pentab:2} for the problem 5 when $P^*$ is equal to $15$.
\end{enumerate}

It must be noticed that in  Tables \ref{pentab:1}, \ref{pentab:2}
the meaning of the column ``Eval." is different in comparison
with Tables \ref{tabpr:1} and II. In  Tables \ref{pentab:1},
\ref{pentab:2} this column shows the total number of evaluations
of the objective function $f(x)$ and {\it all} the constraints.
Thus, it is equal to
\[  (N_v +1)\times N_{iter}, \]
where $N_v$ is the number of constraints and $N_{iter}$ is the
number of iterations for each problem.

In Figures 4 and 5 we show the dynamic diagrams of the search
executed by the IBBA and the PEN for the differentiable problem 7
from \cite{Famularo_Sergeyev_Pugliese}:
\[
\ds \min_{x \in \left[ -3, 2 \right]}  f(x)  =  \ds
\exp\left(-\cos\left(4 x - 3\right)\right)+\frac{1}{250} \left( 4
x - 3 \right)^2-1
 \]
  subject to
\[
\begin{array}{cccl}
& g_1(x) & = & \ds \sin^3(x) \exp(-\sin(3 x))+\frac{1}{2} \le 0 \\[12pt]
& g_2(x) & = & \ds \cos\left(\frac{7}{5} (x+3)
\right)-\sin\left(7 (x+3) \right)+\frac{3}{10} \le 0
\end{array}
\]
The problem has two disjoint feasible subregions shown by two
continuous bold lines and the global optimum $x^*$ is located at
the point $x^*=-0.774575$.

    \begin{center}
    \begin{table}[ht]
\caption{Differentiable problems: comparison between the IBBA and
the~PEN.}
    \label{comptab:1}

    \begin{tabular}{l|rrr|rrr} \hline\noalign{\smallskip}
 & \multicolumn{3}{|c}{Iterations} & \multicolumn{3}{|c}{Evaluations}\\
Problem & PEN & IBBA & Speedup & PEN & IBBA & Speedup \\
\noalign{\smallskip}\hline\noalign{\smallskip}
 1 & $83$ & $23$ & $3.61$ & $166$ & $36$ & $4.61$  \\
 2 & $954$ & $227$ & $4.20$ & $1906$ & $248$ & $7.69$  \\
 3 & $119$ & $62$ & $1.92$ & $238$ & $84$ & $2.83$  \\
 4 & $1762$ & $953$ & $1.85$ & $5286$ & $1459$ & $3.62$ \\
 5 & $765$ & $144$ & $5.31$ & $2295$ & $402$ & $5.71$ \\
 6 & $477$ & $112$ & $4.26$ & $1431$ & $228$ & $6.28$ \\
 7 & $917$ & $294$ & $3.12$ & $2751$ & $794$ & $3.46$ \\
 8 & $821$ & $237$ & $3.46$ & $3284$ & $536$ & $6.13$ \\
 9 & $262$ & $100$ & $2.62$ & $1048$ & $301$ & $3.48$ \\
 10 & $2019$ & $1420$ & $1.42$ & $8076$ & $5344$ & $1.51$ \\
  Average & $817.9$ & $357.2$ & $2.29$ & $2648.1$ & $943.2$ & $2.81$
\\ \noalign{\smallskip}\hline\noalign{\smallskip}
   \end{tabular}
\end{table}
\end{center}

   \begin{center}
    \begin{table}
   \caption{Non-differentiable problems: comparison between the the IBBA and
the~PEN.}
    \label{ncomptab:1}

\begin{tabular}{l|rrr|rrr} \hline\noalign{\smallskip}
 & \multicolumn{3}{|c}{Iterations} & \multicolumn{3}{|c}{Evaluations}\\
Problem & PEN & IBBA & Speedup & PEN & IBBA & Speedup \\
\noalign{\smallskip}\hline\noalign{\smallskip}
 1 &  $247$ & $51$ & $4.84$ & $494$ & $79$ & $6.25$ \\
 2 &  $241$ & $34$ & $7.09$ & $482$ & $50$ & $9.64$ \\
 3 &  $797$ & $190$ & $4.19$ & $1594$ & $209$ & $7.63$ \\
 4 &  $272$ & $235$ & $1.16$ & $819$ & $418$ & $1.96$ \\
 5 &  $671$ & $283$ & $2.37$ & $2013$ & $422$ & $4.77$ \\
 6 &  $909$ & $629$ & $1.45$ & $2727$ & $1839$ & $1.48$ \\
 7 &  $199$ & $120$ & $1.66$ & $597$ & $216$ & $2.76$ \\
 8 &  $365$ & $64$ & $5.70$ & $1460$ & $144$ & $10.14$ \\
 9 &  $1183$ & $334$ & $3.54$ & $4732$ & $1083$ & $4.37$ \\
 10 & $135$ & $75$ & $1.80$ & $540$ & $151$ & $3.58$ \\
 Average & $501.9$ & $201.5$ & $2.49$ & $1545.8$ & $461.1$ & $3.35$ \\
\noalign{\smallskip}\hline\noalign{\smallskip}
   \end{tabular}

\end{table}
   \end{center}

\vspace{-1cm}

The first line (from up to down) of ``+'' located under the graph
of the problem 7 in the upper subplot of Figure~4 represents the
points where the first constraint has not been satisfied (number
of iterations equal to 35). Thus, due to the decision rules of the
IBBA, the second constraint has not been evaluated at these
points. The second line of ``+'' represents the points where the
first constraint has been satisfied but the second constraint has
been not (number of iterations equal to 18). In these points both
constraints have been evaluated but the objective function has
been not. The last line represents the points where both the
constraints have been satisfied (number of evaluations equal to
241). The total number of evaluations is equal to $35+18 \times
2+241 \times 3 = 794$. These evaluations have been executed
during $35+18+241 = 294$ iterations.

The line of ``+'' located under the graph in the upper subplot of
Figure 5 represents the points where the function (\ref{pen}) has
been evaluated. The number of iterations is equal to $917$ and
the number of evaluations is equal to $917 \times 3 = 2757$.

Finally, the infeasibility of the differentiable problem  from
\cite{Famularo_Sergeyev_Pugliese} has been determined by the IBBA
in $38$ iterations consisting of $9$ evaluations of the first
constraint and $29$ evaluations of the first and second
constraints (i.e., $67$ evaluations in total). The infeasibility
of the non-differentiable problem   from
\cite{Famularo_Sergeyev_Pugliese} has been determined by the IBBA
in $98$ iterations consisting of $93$ evaluations of the first
constraint and $5$ evaluations of the first and second
constraints (i.e., $103$ evaluations in total). Naturally, the
objective functions were not evaluated at all in both cases. Note
that experiments for infeasible problems have not been executed
with the PEN because the penalty approach does not allow to the
user to determine infeasibility of problems.

\section{Concluding remarks}
Lipschitz univariate constrained global optimization problems
where both the objective function and  constraints  can be
multiextremal  have been considered in this paper. The
constrained problem has been reduced to a discontinuous
unconstrained problem by the index scheme. A Branch-and-Bound
method for solving the reduced problem has been  proposed.
Convergence conditions of the new method have been established.

The new algorithm works without usage of derivatives. It either
determines the infeasibility of the original problems or finds
upper and lower bounds of the global solution. Note that it is
able to work with problems where the objective function and/or
constraints are not defined over the whole search region. It does
not evaluate all the constraints during every iteration. The
introduction of additional variables and/or parameters is not
required.

Extensive numerical results show quite a satisfactory performance
of the new technique. The behaviour of the Index Branch-and-Bound
method was compared to the method of Pijavskii combined with a
penalty approach. This algorithm has been chosen for comparison
because it used the same information about the problem as the
IBBA -- the Lipschitz constants for the objective function and
constraints.

A priori the penalty approach combined with the method of
Pijavskii seemed to be more attractive because it dealt only with
one function. In the facts however, the evaluation of this
function requires the evaluation of $m+1$ initial functions. The
second disadvantage of the penalty approach is that it requires
an accurate tuning of the penalty coefficient in contrast to the
IBBA which works without any additional parameter. Finally, when
 the penalty approach is used and a constraint $g(x)$ is
defined only over a subregion $[c,d]$ of the search region
$[a,b]$, the problem of extending  $g(x)$ to the whole region
$[a,b]$ arises. In contrast, the IBBA does not have this
difficulty because every constraint (and the objective function)
is evaluated only within its region of definition.

\begin{figure}[ht]
  \begin{center}
    \caption{Optimization of the differentiable problem 7 by the IBBA.  \label{figpr:7}}
    \epsfig{ figure = 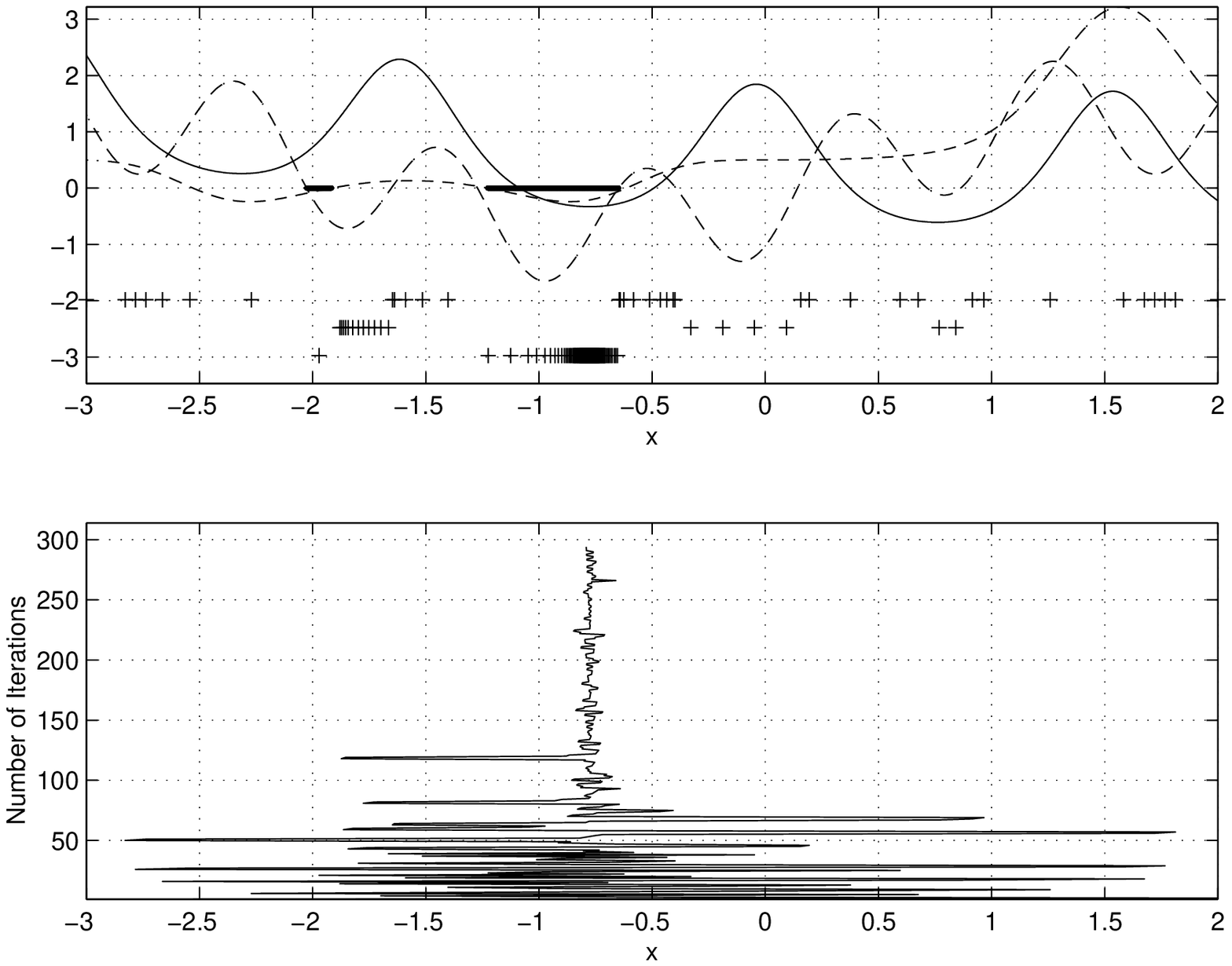, width = 4.5in, height = 3.375in, silent = yes }
  \end{center}
\end{figure}
\begin{figure}[h]
  \begin{center}
    \caption{Optimization of the differentiable problem 7  by the PEN.  \label{figpr:Pijavskii}}
    \epsfig{ figure = 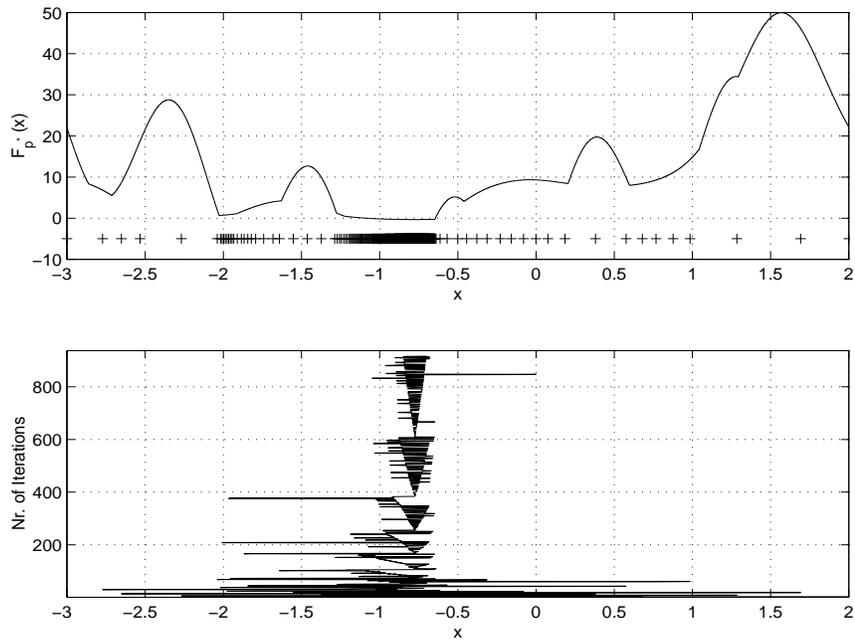, width = 4.5in, height = 3.375in, silent = yes }
  \end{center}
\end{figure}

{}

\end{article}
\end{document}

%% file: fig1.eepic.tex
\setlength{\unitlength}{0.00041667in}
\begingroup\makeatletter\ifx\SetFigFont\undefined
\def\x#1#2#3#4#5#6#7\relax{\def\x{#1#2#3#4#5#6}}%
\expandafter\x\fmtname xxxxxx\relax \def\y{splain}%
\ifx\x\y   
\gdef\SetFigFont#1#2#3{%
  \ifnum #1<17\tiny\else \ifnum #1<20\small\else
  \ifnum #1<24\normalsize\else \ifnum #1<29\large\else
  \ifnum #1<34\Large\else \ifnum #1<41\LARGE\else
     \huge\fi\fi\fi\fi\fi\fi
  \csname #3\endcsname}%
\else
\gdef\SetFigFont#1#2#3{\begingroup
  \count@#1\relax \ifnum 25<\count@\count@25\fi
  \def\x{\endgroup\@setsize\SetFigFont{#2pt}}%
  \expandafter\x
    \csname \romannumeral\the\count@ pt\expandafter\endcsname
    \csname @\romannumeral\the\count@ pt\endcsname
  \csname #3\endcsname}%
\fi
\fi\endgroup
{\renewcommand{\dashlinestretch}{30}
\begin{picture}(11424,9789)(0,-10)
\put(6162,1992){\makebox(0,0)[lb]{\smash{{{\SetFigFont{6}{7.2}{rm}$Q_3$}}}}}
\put(1662,1992){\makebox(0,0)[lb]{\smash{{{\SetFigFont{6}{7.2}{rm}$Q_3$}}}}}
\put(6162,1092){\makebox(0,0)[lb]{\smash{{{\SetFigFont{6}{7.2}{rm}$Q_2$}}}}}
\put(1512,1092){\makebox(0,0)[lb]{\smash{{{\SetFigFont{6}{7.2}{rm}$Q_2$}}}}}
\put(10512,1092){\makebox(0,0)[lb]{\smash{{{\SetFigFont{6}{7.2}{rm}$Q_2$}}}}}
\put(11187,4287){\makebox(0,0)[lb]{\smash{{{\SetFigFont{6}{7.2}{rm}b}}}}}
\put(1062,4287){\makebox(0,0)[lb]{\smash{{{\SetFigFont{6}{7.2}{rm}a}}}}}
\put(7587,4287){\makebox(0,0)[lb]{\smash{{{\SetFigFont{6}{7.2}{rm}c}}}}}
\put(10062,4287){\makebox(0,0)[lb]{\smash{{{\SetFigFont{6}{7.2}{rm}d}}}}}
\put(1212,6912){\blacken\ellipse{94}{94}}
\put(1212,6912){\ellipse{94}{94}}
\put(1512,7512){\blacken\ellipse{94}{94}}
\put(1512,7512){\ellipse{94}{94}}
\put(2112,8097){\blacken\ellipse{94}{94}}
\put(2112,8097){\ellipse{94}{94}}
\put(5412,8412){\blacken\ellipse{94}{94}}
\put(5412,8412){\ellipse{94}{94}}
\put(7212,9612){\blacken\ellipse{84}{84}}
\put(7212,9612){\ellipse{84}{84}}
\put(7512,5412){\blacken\ellipse{94}{94}}
\put(7512,5412){\ellipse{94}{94}}
\put(10212,6912){\blacken\ellipse{94}{94}}
\put(10212,6912){\ellipse{94}{94}}
\put(11112,7512){\blacken\ellipse{84}{84}}
\put(11112,7512){\ellipse{84}{84}}
\path(12,4512)(11412,4512)
\blacken\path(11292.000,4482.000)(11412.000,4512.000)(11292.000,4542.000)(11292.000,4482.000)
\dashline{60.000}(1512,7512)(1512,1812)
\dashline{60.000}(1212,6912)(1212,12)
\dashline{60.000}(7212,1812)(7212,9612)
\dashline{60.000}(11112,7512)(11112,12)
\dashline{60.000}(2112,8112)(2112,912)
\dashline{60.000}(7512,5412)(7512,912)
\dashline{60.000}(5412,8412)(5412,912)
\dashline{60.000}(10212,912)(10212,6912)
\path(1512,1812)(2112,1812)
\blacken\path(1992.000,1782.000)(2112.000,1812.000)(1992.000,1842.000)(1992.000,1782.000)
\path(2112,1812)(1512,1812)
\blacken\path(1632.000,1842.000)(1512.000,1812.000)(1632.000,1782.000)(1632.000,1842.000)
\path(1212,912)(2112,912)
\blacken\path(1992.000,882.000)(2112.000,912.000)(1992.000,942.000)(1992.000,882.000)
\path(2112,912)(1212,912)
\blacken\path(1332.000,942.000)(1212.000,912.000)(1332.000,882.000)(1332.000,942.000)
\path(5412,1812)(7212,1812)
\blacken\path(7092.000,1782.000)(7212.000,1812.000)(7092.000,1842.000)(7092.000,1782.000)
\path(7212,1812)(5412,1812)
\blacken\path(5532.000,1842.000)(5412.000,1812.000)(5532.000,1782.000)(5532.000,1842.000)
\path(5412,912)(7512,912)
\blacken\path(7392.000,882.000)(7512.000,912.000)(7392.000,942.000)(7392.000,882.000)
\path(7512,912)(5412,912)
\blacken\path(5532.000,942.000)(5412.000,912.000)(5532.000,882.000)(5532.000,942.000)
\path(10212,912)(11112,912)
\blacken\path(10992.000,882.000)(11112.000,912.000)(10992.000,942.000)(10992.000,882.000)
\path(11112,912)(10212,912)
\blacken\path(10332.000,942.000)(10212.000,912.000)(10332.000,882.000)(10332.000,942.000)
\path(1212,12)(11112,12)
\blacken\path(10992.000,-18.000)(11112.000,12.000)(10992.000,42.000)(10992.000,-18.000)
\path(11112,12)(1212,12)
\blacken\path(1332.000,42.000)(1212.000,12.000)(1332.000,-18.000)(1332.000,42.000)
\dashline{60.000}(6387,6987)(6387,4512)
\path(612,3912)(612,9762)
\path(612,3912)(612,9762)
\blacken\path(642.000,9642.000)(612.000,9762.000)(582.000,9642.000)(642.000,9642.000)
\thicklines
\path(10212,6912)(10212,6913)(10214,6917)
	(10217,6928)(10224,6947)(10233,6974)
	(10246,7009)(10260,7049)(10276,7092)
	(10293,7135)(10310,7176)(10326,7214)
	(10342,7246)(10356,7273)(10370,7294)
	(10383,7310)(10394,7319)(10406,7322)
	(10416,7320)(10427,7312)(10437,7299)
	(10443,7289)(10450,7276)(10456,7261)
	(10463,7244)(10470,7224)(10477,7201)
	(10484,7175)(10491,7146)(10499,7114)
	(10507,7079)(10516,7040)(10525,6997)
	(10534,6952)(10544,6903)(10554,6850)
	(10564,6796)(10574,6739)(10585,6681)
	(10596,6622)(10606,6564)(10616,6508)
	(10625,6455)(10633,6407)(10641,6363)
	(10647,6327)(10652,6296)(10656,6273)
	(10659,6256)(10661,6245)(10662,6240)(10662,6237)
\path(1512,7512)(1515,7509)(1522,7502)
	(1533,7491)(1549,7476)(1567,7458)
	(1587,7438)(1608,7418)(1627,7399)
	(1646,7382)(1662,7367)(1677,7353)
	(1690,7341)(1702,7331)(1714,7321)
	(1725,7312)(1732,7306)(1739,7301)
	(1746,7296)(1753,7291)(1760,7287)
	(1767,7283)(1774,7280)(1782,7277)
	(1789,7276)(1797,7275)(1804,7275)
	(1812,7275)(1820,7277)(1827,7280)
	(1835,7285)(1842,7290)(1850,7296)
	(1857,7304)(1864,7312)(1871,7322)
	(1878,7333)(1885,7346)(1892,7359)
	(1900,7374)(1906,7390)(1913,7406)
	(1920,7425)(1927,7445)(1935,7468)
	(1943,7493)(1952,7521)(1961,7551)
	(1971,7585)(1981,7621)(1992,7661)
	(2004,7703)(2016,7748)(2029,7794)
	(2041,7842)(2054,7889)(2066,7934)
	(2077,7976)(2087,8014)(2095,8046)
	(2101,8071)(2106,8090)(2109,8102)
	(2111,8109)(2112,8112)
\path(5412,8412)(5414,8410)(5419,8405)
	(5427,8397)(5439,8384)(5455,8367)
	(5475,8347)(5497,8323)(5521,8298)
	(5545,8272)(5569,8246)(5592,8221)
	(5614,8197)(5635,8174)(5653,8152)
	(5671,8132)(5686,8112)(5701,8094)
	(5714,8076)(5727,8059)(5738,8042)
	(5750,8024)(5763,8003)(5775,7982)
	(5788,7960)(5800,7938)(5812,7916)
	(5824,7894)(5836,7873)(5847,7852)
	(5859,7833)(5870,7815)(5881,7798)
	(5891,7784)(5901,7772)(5911,7763)
	(5921,7756)(5931,7751)(5940,7749)
	(5950,7749)(5959,7753)(5969,7758)
	(5980,7766)(5990,7775)(6002,7786)
	(6014,7798)(6026,7810)(6038,7822)
	(6051,7834)(6064,7844)(6077,7853)
	(6090,7860)(6103,7865)(6115,7866)
	(6127,7864)(6139,7859)(6151,7850)
	(6162,7837)(6169,7826)(6177,7813)
	(6184,7798)(6192,7780)(6200,7760)
	(6208,7737)(6216,7710)(6225,7680)
	(6235,7647)(6245,7610)(6255,7570)
	(6266,7526)(6277,7480)(6289,7430)
	(6301,7378)(6313,7325)(6325,7272)
	(6336,7220)(6347,7171)(6357,7127)
	(6366,7087)(6373,7054)(6378,7028)
	(6382,7009)(6385,6997)(6386,6990)(6387,6987)
\path(6387,6987)(6387,6990)(6388,6996)
	(6390,7008)(6392,7026)(6395,7051)
	(6399,7083)(6405,7121)(6411,7166)
	(6417,7216)(6425,7269)(6432,7325)
	(6440,7381)(6449,7438)(6457,7493)
	(6465,7546)(6473,7597)(6480,7645)
	(6488,7689)(6495,7730)(6502,7768)
	(6509,7802)(6516,7833)(6522,7861)
	(6529,7886)(6535,7908)(6542,7928)
	(6549,7946)(6555,7961)(6562,7974)
	(6575,7995)(6589,8011)(6604,8022)
	(6619,8029)(6635,8032)(6652,8033)
	(6669,8031)(6687,8028)(6705,8024)
	(6722,8021)(6739,8017)(6755,8016)
	(6770,8016)(6785,8020)(6799,8026)
	(6812,8037)(6820,8046)(6828,8056)
	(6835,8069)(6843,8083)(6850,8100)
	(6858,8118)(6865,8138)(6873,8160)
	(6880,8183)(6888,8209)(6895,8235)
	(6903,8264)(6910,8293)(6918,8324)
	(6925,8355)(6932,8387)(6939,8420)
	(6946,8453)(6953,8486)(6960,8518)
	(6967,8551)(6973,8584)(6980,8616)
	(6986,8648)(6993,8680)(7000,8712)
	(7005,8740)(7011,8767)(7017,8795)
	(7023,8823)(7030,8853)(7037,8883)
	(7044,8914)(7051,8946)(7059,8980)
	(7067,9016)(7076,9053)(7085,9092)
	(7095,9133)(7106,9176)(7116,9220)
	(7127,9265)(7138,9310)(7149,9354)
	(7159,9398)(7169,9439)(7179,9477)
	(7187,9511)(7194,9540)(7200,9565)
	(7205,9583)(7208,9597)(7210,9605)
	(7211,9610)(7212,9612)
\path(10662,6237)(10662,6239)(10663,6245)
	(10664,6255)(10665,6270)(10667,6292)
	(10670,6320)(10673,6356)(10677,6398)
	(10681,6447)(10687,6502)(10692,6561)
	(10698,6623)(10705,6689)(10711,6755)
	(10718,6823)(10725,6889)(10732,6955)
	(10738,7018)(10745,7079)(10752,7137)
	(10758,7192)(10764,7244)(10770,7293)
	(10776,7338)(10782,7380)(10788,7419)
	(10793,7454)(10799,7487)(10804,7517)
	(10809,7545)(10815,7570)(10820,7593)
	(10826,7614)(10831,7632)(10837,7649)
	(10845,7670)(10853,7688)(10862,7702)
	(10871,7714)(10881,7721)(10891,7726)
	(10902,7728)(10914,7726)(10926,7721)
	(10940,7713)(10954,7703)(10970,7690)
	(10986,7674)(11002,7656)(11019,7637)
	(11035,7617)(11051,7597)(11065,7578)
	(11078,7560)(11089,7545)(11098,7533)
	(11104,7523)(11109,7517)(11111,7514)(11112,7512)
\path(1212,6912)(1212,6910)(1213,6906)
	(1213,6897)(1214,6884)(1216,6866)
	(1218,6841)(1221,6810)(1224,6773)
	(1228,6730)(1233,6682)(1237,6628)
	(1243,6571)(1248,6511)(1254,6448)
	(1260,6384)(1266,6319)(1271,6254)
	(1277,6190)(1283,6128)(1289,6067)
	(1294,6009)(1300,5953)(1305,5899)
	(1310,5847)(1315,5799)(1319,5752)
	(1324,5708)(1328,5666)(1332,5626)
	(1336,5588)(1340,5552)(1344,5518)
	(1347,5484)(1351,5452)(1355,5421)
	(1358,5391)(1362,5362)(1367,5323)
	(1372,5286)(1377,5249)(1382,5213)
	(1387,5177)(1393,5142)(1398,5106)
	(1404,5070)(1410,5033)(1417,4996)
	(1424,4957)(1431,4918)(1439,4877)
	(1447,4836)(1455,4795)(1463,4754)
	(1471,4714)(1478,4676)(1486,4640)
	(1492,4608)(1498,4580)(1502,4558)
	(1506,4540)(1509,4527)(1511,4519)
	(1512,4514)(1512,4512)
\path(1512,4512)(1513,4510)(1514,4505)
	(1517,4496)(1522,4482)(1529,4463)
	(1538,4438)(1548,4408)(1561,4372)
	(1575,4331)(1591,4286)(1608,4239)
	(1625,4190)(1643,4139)(1662,4089)
	(1680,4039)(1698,3990)(1715,3943)
	(1732,3898)(1748,3856)(1764,3815)
	(1779,3777)(1794,3740)(1808,3706)
	(1822,3674)(1835,3644)(1848,3615)
	(1861,3587)(1874,3561)(1886,3535)
	(1899,3511)(1912,3487)(1928,3458)
	(1944,3430)(1961,3403)(1978,3375)
	(1995,3349)(2013,3322)(2032,3297)
	(2051,3271)(2070,3246)(2090,3222)
	(2109,3198)(2129,3176)(2150,3154)
	(2170,3133)(2190,3113)(2209,3095)
	(2229,3077)(2248,3061)(2267,3046)
	(2286,3033)(2304,3020)(2321,3009)
	(2338,2999)(2355,2990)(2371,2982)
	(2387,2974)(2406,2967)(2424,2961)
	(2443,2956)(2461,2953)(2480,2950)
	(2499,2949)(2518,2949)(2537,2951)
	(2555,2954)(2574,2958)(2593,2963)
	(2611,2970)(2628,2977)(2646,2986)
	(2662,2996)(2678,3006)(2694,3018)
	(2709,3030)(2723,3044)(2736,3057)
	(2749,3072)(2762,3087)(2773,3102)
	(2785,3117)(2796,3134)(2807,3152)
	(2818,3171)(2830,3193)(2842,3216)
	(2854,3241)(2867,3268)(2881,3298)
	(2895,3330)(2910,3363)(2925,3398)
	(2939,3433)(2954,3468)(2967,3501)
	(2980,3531)(2990,3557)(2998,3578)
	(3005,3593)(3009,3604)(3011,3609)(3012,3612)
\path(3012,3612)(3015,3613)(3021,3615)
	(3032,3618)(3048,3622)(3069,3627)
	(3093,3632)(3121,3637)(3149,3642)
	(3178,3645)(3205,3647)(3232,3648)
	(3257,3646)(3280,3643)(3302,3637)
	(3323,3630)(3343,3620)(3362,3607)
	(3381,3592)(3400,3574)(3413,3560)
	(3427,3545)(3441,3528)(3455,3510)
	(3469,3491)(3484,3470)(3499,3448)
	(3515,3424)(3530,3400)(3546,3375)
	(3563,3348)(3579,3321)(3595,3294)
	(3612,3266)(3629,3238)(3645,3209)
	(3661,3181)(3678,3153)(3694,3125)
	(3709,3098)(3725,3072)(3740,3046)
	(3755,3022)(3769,2998)(3783,2975)
	(3797,2953)(3811,2932)(3825,2912)
	(3842,2888)(3859,2865)(3876,2843)
	(3894,2822)(3912,2802)(3932,2782)
	(3952,2762)(3974,2742)(3997,2722)
	(4021,2702)(4046,2682)(4072,2662)
	(4098,2642)(4123,2624)(4147,2607)
	(4167,2593)(4184,2581)(4197,2572)
	(4205,2567)(4210,2563)(4212,2562)
\path(4212,2562)(4214,2561)(4220,2558)
	(4229,2553)(4244,2545)(4264,2534)
	(4290,2520)(4321,2504)(4356,2485)
	(4396,2464)(4437,2442)(4480,2420)
	(4524,2397)(4567,2374)(4609,2353)
	(4650,2332)(4689,2312)(4726,2293)
	(4761,2276)(4794,2260)(4826,2244)
	(4855,2230)(4884,2217)(4911,2205)
	(4937,2193)(4962,2182)(4987,2172)
	(5012,2162)(5038,2152)(5065,2142)
	(5091,2132)(5118,2123)(5145,2114)
	(5173,2105)(5201,2097)(5229,2089)
	(5258,2082)(5287,2075)(5317,2069)
	(5346,2064)(5376,2060)(5406,2056)
	(5436,2053)(5466,2050)(5495,2049)
	(5525,2049)(5554,2049)(5583,2050)
	(5611,2052)(5639,2055)(5667,2059)
	(5694,2063)(5722,2068)(5750,2074)
	(5774,2081)(5798,2087)(5823,2095)
	(5848,2103)(5874,2112)(5901,2122)
	(5928,2133)(5955,2145)(5983,2158)
	(6011,2172)(6040,2187)(6068,2202)
	(6097,2219)(6126,2236)(6155,2254)
	(6184,2273)(6212,2292)(6240,2313)
	(6268,2333)(6295,2354)(6321,2376)
	(6347,2398)(6372,2421)(6396,2443)
	(6419,2467)(6442,2490)(6463,2514)
	(6484,2538)(6505,2562)(6525,2587)
	(6543,2611)(6560,2635)(6578,2661)
	(6595,2687)(6612,2713)(6628,2741)
	(6645,2769)(6661,2799)(6677,2829)
	(6693,2860)(6709,2891)(6724,2924)
	(6740,2957)(6755,2990)(6770,3024)
	(6784,3058)(6798,3092)(6812,3127)
	(6825,3161)(6838,3196)(6850,3230)
	(6862,3263)(6874,3297)(6885,3330)
	(6896,3362)(6906,3394)(6916,3425)
	(6926,3456)(6935,3486)(6944,3516)
	(6953,3545)(6962,3574)(6971,3605)
	(6980,3636)(6989,3667)(6998,3698)
	(7007,3729)(7016,3762)(7025,3795)
	(7034,3829)(7044,3865)(7054,3902)
	(7064,3941)(7075,3982)(7086,4024)
	(7098,4068)(7110,4113)(7122,4159)
	(7133,4205)(7145,4251)(7157,4295)
	(7167,4337)(7177,4376)(7186,4410)
	(7194,4440)(7200,4464)(7205,4483)
	(7208,4496)(7210,4505)(7211,4510)(7212,4512)
\path(7212,4512)(7213,4514)(7214,4518)
	(7217,4526)(7221,4539)(7227,4557)
	(7235,4581)(7245,4611)(7257,4647)
	(7271,4689)(7287,4736)(7304,4788)
	(7323,4844)(7342,4902)(7362,4962)
	(7382,5022)(7401,5080)(7420,5136)
	(7437,5188)(7453,5235)(7467,5277)
	(7479,5313)(7489,5343)(7497,5367)
	(7503,5385)(7507,5398)(7510,5406)
	(7511,5410)(7512,5412)
\path(2112,4512)(2114,4513)(2118,4514)
	(2125,4516)(2137,4519)(2153,4523)
	(2173,4529)(2198,4536)(2226,4545)
	(2258,4554)(2292,4565)(2328,4576)
	(2365,4587)(2402,4599)(2439,4612)
	(2475,4624)(2509,4637)(2542,4650)
	(2573,4663)(2603,4676)(2631,4690)
	(2659,4705)(2685,4720)(2712,4737)
	(2733,4751)(2755,4766)(2777,4782)
	(2799,4798)(2821,4815)(2844,4833)
	(2867,4851)(2891,4869)(2914,4888)
	(2938,4908)(2962,4927)(2987,4947)
	(3011,4968)(3036,4988)(3060,5009)
	(3084,5030)(3109,5051)(3133,5072)
	(3157,5093)(3180,5114)(3203,5136)
	(3226,5157)(3248,5179)(3270,5201)
	(3291,5223)(3312,5245)(3331,5268)
	(3351,5290)(3369,5313)(3387,5337)
	(3404,5361)(3421,5386)(3437,5412)
	(3453,5439)(3470,5469)(3486,5499)
	(3502,5532)(3519,5567)(3536,5603)
	(3553,5641)(3570,5681)(3587,5722)
	(3605,5764)(3622,5806)(3640,5848)
	(3656,5890)(3672,5930)(3687,5969)
	(3701,6004)(3714,6037)(3725,6066)
	(3735,6092)(3743,6113)(3750,6130)
	(3755,6142)(3758,6152)(3760,6157)
	(3761,6161)(3762,6162)
\path(1212,3312)(1213,3313)(1214,3317)
	(1217,3323)(1222,3332)(1229,3346)
	(1238,3363)(1249,3385)(1262,3412)
	(1277,3442)(1295,3476)(1313,3513)
	(1334,3553)(1355,3594)(1377,3637)
	(1399,3680)(1422,3724)(1444,3766)
	(1466,3808)(1488,3848)(1508,3886)
	(1528,3923)(1548,3957)(1566,3989)
	(1583,4018)(1600,4046)(1616,4071)
	(1632,4095)(1647,4116)(1662,4137)
	(1681,4161)(1700,4184)(1720,4206)
	(1741,4228)(1763,4249)(1786,4270)
	(1810,4292)(1836,4313)(1863,4334)
	(1891,4355)(1919,4376)(1947,4397)
	(1974,4417)(2001,4435)(2025,4452)
	(2047,4467)(2065,4480)(2081,4491)
	(2093,4499)(2102,4505)(2107,4509)
	(2111,4511)(2112,4512)
\path(3762,6162)(3762,6161)(3764,6158)
	(3766,6152)(3769,6144)(3773,6132)
	(3780,6116)(3787,6096)(3797,6072)
	(3808,6044)(3821,6011)(3835,5975)
	(3851,5935)(3868,5893)(3886,5848)
	(3905,5802)(3925,5754)(3945,5706)
	(3966,5658)(3987,5610)(4008,5563)
	(4029,5518)(4050,5473)(4071,5431)
	(4092,5390)(4113,5352)(4133,5315)
	(4154,5281)(4175,5249)(4196,5218)
	(4218,5190)(4240,5162)(4263,5137)
	(4287,5112)(4312,5088)(4338,5065)
	(4366,5042)(4396,5019)(4427,4997)
	(4461,4974)(4498,4951)(4536,4928)
	(4577,4905)(4621,4881)(4666,4857)
	(4714,4833)(4763,4808)(4814,4783)
	(4865,4758)(4917,4734)(4970,4709)
	(5021,4686)(5071,4663)(5120,4641)
	(5166,4620)(5209,4601)(5248,4584)
	(5283,4568)(5314,4555)(5340,4543)
	(5362,4534)(5379,4526)(5392,4520)
	(5402,4516)(5407,4514)(5411,4513)(5412,4512)
\path(5412,4512)(5414,4511)(5419,4510)
	(5427,4507)(5440,4503)(5458,4497)
	(5482,4489)(5512,4479)(5546,4468)
	(5585,4456)(5628,4443)(5674,4428)
	(5722,4414)(5771,4399)(5821,4384)
	(5870,4370)(5918,4357)(5964,4344)
	(6009,4332)(6053,4321)(6094,4311)
	(6134,4302)(6172,4294)(6208,4287)
	(6243,4281)(6276,4275)(6309,4271)
	(6340,4268)(6371,4265)(6402,4263)
	(6432,4262)(6462,4262)(6492,4262)
	(6522,4263)(6553,4265)(6584,4268)
	(6615,4271)(6648,4275)(6681,4281)
	(6716,4287)(6752,4294)(6790,4302)
	(6830,4311)(6871,4321)(6915,4332)
	(6960,4344)(7006,4357)(7054,4370)
	(7103,4384)(7153,4399)(7202,4414)
	(7250,4428)(7296,4443)(7339,4456)
	(7378,4468)(7412,4479)(7442,4489)
	(7466,4497)(7484,4503)(7497,4507)
	(7505,4510)(7510,4511)(7512,4512)
\path(7512,4512)(7513,4514)(7515,4517)
	(7519,4524)(7526,4534)(7535,4549)
	(7547,4568)(7561,4592)(7579,4620)
	(7599,4652)(7621,4688)(7645,4727)
	(7670,4768)(7697,4811)(7724,4854)
	(7751,4898)(7778,4941)(7804,4984)
	(7830,5026)(7855,5066)(7880,5105)
	(7903,5143)(7926,5179)(7948,5213)
	(7968,5247)(7989,5279)(8008,5310)
	(8027,5340)(8046,5370)(8064,5399)
	(8083,5427)(8101,5456)(8119,5484)
	(8137,5512)(8156,5542)(8176,5573)
	(8196,5603)(8216,5634)(8237,5665)
	(8258,5696)(8279,5728)(8301,5759)
	(8324,5791)(8347,5823)(8370,5854)
	(8393,5885)(8417,5915)(8442,5945)
	(8466,5974)(8490,6003)(8515,6030)
	(8540,6056)(8564,6081)(8589,6104)
	(8613,6126)(8637,6146)(8661,6165)
	(8684,6181)(8707,6196)(8730,6209)
	(8753,6221)(8775,6230)(8797,6238)
	(8819,6244)(8840,6247)(8862,6249)
	(8885,6250)(8908,6248)(8931,6243)
	(8955,6237)(8979,6229)(9003,6218)
	(9027,6206)(9052,6191)(9077,6175)
	(9103,6156)(9129,6135)(9155,6113)
	(9181,6089)(9207,6064)(9234,6037)
	(9260,6009)(9286,5979)(9311,5949)
	(9337,5918)(9362,5886)(9387,5854)
	(9411,5822)(9434,5789)(9457,5756)
	(9480,5724)(9502,5691)(9524,5659)
	(9545,5626)(9566,5594)(9587,5562)
	(9605,5534)(9623,5505)(9641,5477)
	(9660,5448)(9678,5418)(9697,5388)
	(9716,5357)(9735,5325)(9756,5291)
	(9776,5256)(9798,5220)(9821,5182)
	(9844,5143)(9869,5102)(9894,5059)
	(9920,5015)(9946,4970)(9973,4924)
	(10000,4877)(10027,4831)(10054,4785)
	(10079,4742)(10103,4700)(10125,4662)
	(10145,4627)(10163,4597)(10177,4572)
	(10189,4551)(10198,4536)(10205,4525)
	(10209,4518)(10211,4514)(10212,4512)
\path(10212,4512)(10213,4509)(10216,4504)
	(10222,4493)(10230,4477)(10241,4456)
	(10254,4430)(10270,4400)(10287,4366)
	(10305,4331)(10324,4295)(10343,4259)
	(10361,4225)(10378,4192)(10395,4161)
	(10411,4132)(10426,4105)(10439,4081)
	(10453,4058)(10465,4036)(10477,4016)
	(10489,3998)(10501,3979)(10512,3962)
	(10526,3942)(10540,3922)(10554,3902)
	(10568,3883)(10583,3864)(10598,3845)
	(10614,3826)(10630,3808)(10646,3790)
	(10662,3772)(10678,3756)(10694,3739)
	(10710,3724)(10726,3709)(10741,3695)
	(10756,3682)(10770,3670)(10784,3658)
	(10798,3648)(10812,3637)(10827,3626)
	(10843,3615)(10859,3604)(10876,3593)
	(10894,3582)(10913,3570)(10934,3558)
	(10957,3545)(10981,3532)(11006,3518)
	(11031,3505)(11054,3492)(11075,3481)
	(11091,3473)(11102,3467)(11109,3464)(11112,3462)
\put(6162,162){\makebox(0,0)[lb]{\smash{{{\SetFigFont{6}{7.2}{rm}$Q_1$}}}}}
\put(1512,7812){\makebox(0,0)[lb]{\smash{{{\SetFigFont{6}{7.2}{rm}$f(x)$}}}}}
\put(5862,8337){\makebox(0,0)[lb]{\smash{{{\SetFigFont{6}{7.2}{rm}$f(x)$}}}}}
\put(8637,5487){\makebox(0,0)[lb]{\smash{{{\SetFigFont{6}{7.2}{rm}$g_1(x)$}}}}}
\put(2937,3987){\makebox(0,0)[lb]{\smash{{{\SetFigFont{6}{7.2}{rm}$g_2(x)$}}}}}
\put(10287,7887){\makebox(0,0)[lb]{\smash{{{\SetFigFont{6}{7.2}{rm}$g_2(x)$}}}}}
\put(6462,4587){\makebox(0,0)[lb]{\smash{{{\SetFigFont{6}{7.2}{rm}$x^*$}}}}}
\end{picture}
}

%% file: fig2.eepic.tex
\setlength{\unitlength}{0.00041667in}
\begingroup\makeatletter\ifx\SetFigFont\undefined
\def\x#1#2#3#4#5#6#7\relax{\def\x{#1#2#3#4#5#6}}%
\expandafter\x\fmtname xxxxxx\relax \def\y{splain}%
\ifx\x\y   
\gdef\SetFigFont#1#2#3{%
  \ifnum #1<17\tiny\else \ifnum #1<20\small\else
  \ifnum #1<24\normalsize\else \ifnum #1<29\large\else
  \ifnum #1<34\Large\else \ifnum #1<41\LARGE\else
     \huge\fi\fi\fi\fi\fi\fi
  \csname #3\endcsname}%
\else
\gdef\SetFigFont#1#2#3{\begingroup
  \count@#1\relax \ifnum 25<\count@\count@25\fi
  \def\x{\endgroup\@setsize\SetFigFont{#2pt}}%
  \expandafter\x
    \csname \romannumeral\the\count@ pt\expandafter\endcsname
    \csname @\romannumeral\the\count@ pt\endcsname
  \csname #3\endcsname}%
\fi
\fi\endgroup
{\renewcommand{\dashlinestretch}{30}
\begin{picture}(11424,3939)(0,-10)
\put(1212,312){\makebox(0,0)[lb]{\smash{{{\SetFigFont{6}{7.2}{rm}a}}}}}
\put(11112,312){\makebox(0,0)[lb]{\smash{{{\SetFigFont{6}{7.2}{rm}b}}}}}
\put(6312,327){\makebox(0,0)[lb]{\smash{{{\SetFigFont{6}{7.2}{rm}$x^*$}}}}}
\put(1512,1122){\blacken\ellipse{60}{60}}
\put(1512,1122){\ellipse{60}{60}}
\put(2112,1707){\blacken\ellipse{94}{94}}
\put(2112,1707){\ellipse{94}{94}}
\put(1527,1107){\blacken\ellipse{90}{90}}
\put(1527,1107){\ellipse{90}{90}}
\put(5412,2007){\blacken\ellipse{90}{90}}
\put(5412,2007){\ellipse{90}{90}}
\put(7212,3207){\blacken\ellipse{90}{90}}
\put(7212,3207){\ellipse{90}{90}}
\put(7512,1497){\blacken\ellipse{94}{94}}
\put(7512,1497){\ellipse{94}{94}}
\put(10212,3012){\blacken\ellipse{90}{90}}
\put(10212,3012){\ellipse{90}{90}}
\put(11112,3612){\blacken\ellipse{90}{90}}
\put(11112,3612){\ellipse{90}{90}}
\put(1212,3012){\blacken\ellipse{90}{90}}
\put(1212,3012){\ellipse{90}{90}}
\put(6380,612){\blacken\ellipse{90}{90}}
\put(6380,612){\ellipse{90}{90}}
\path(12,612)(11412,612)
\blacken\path(11292.000,582.000)(11412.000,612.000)(11292.000,642.000)(11292.000,582.000)
\dashline{60.000}(1212,3012)(1212,612)
\dashline{60.000}(1512,1062)(1512,612)
\dashline{60.000}(2112,1737)(2112,612)
\dashline{60.000}(5412,2037)(5412,612)
\dashline{60.000}(7212,3237)(7212,612)
\dashline{60.000}(10212,3012)(10212,612)
\dashline{60.000}(11112,3612)(11112,612)
\path(612,12)(612,3912)
\path(612,12)(612,3912)
\blacken\path(642.000,3792.000)(612.000,3912.000)(582.000,3792.000)(642.000,3792.000)
\dashline{60.000}(7512,1512)(7512,612)
\thicklines
\path(1212,3012)(1212,3010)(1213,3006)
    (1213,2997)(1214,2984)(1216,2966)
    (1218,2941)(1221,2910)(1224,2873)
    (1228,2830)(1233,2782)(1237,2728)
    (1243,2671)(1248,2611)(1254,2548)
    (1260,2484)(1266,2419)(1271,2354)
    (1277,2290)(1283,2228)(1289,2167)
    (1294,2109)(1300,2053)(1305,1999)
    (1310,1947)(1315,1899)(1319,1852)
    (1324,1808)(1328,1766)(1332,1726)
    (1336,1688)(1340,1652)(1344,1618)
    (1347,1584)(1351,1552)(1355,1521)
    (1358,1491)(1362,1462)(1367,1423)
    (1372,1386)(1377,1349)(1382,1313)
    (1387,1277)(1393,1242)(1398,1206)
    (1404,1170)(1410,1133)(1417,1096)
    (1424,1057)(1431,1018)(1439,977)
    (1447,936)(1455,895)(1463,854)
    (1471,814)(1478,776)(1486,740)
    (1492,708)(1498,680)(1502,658)
    (1506,640)(1509,627)(1511,619)
    (1512,614)(1512,612)
\path(1512,1137)(1515,1134)(1522,1127)
    (1533,1116)(1549,1101)(1567,1083)
    (1587,1063)(1608,1043)(1627,1024)
    (1646,1007)(1662,992)(1677,978)
    (1690,966)(1702,956)(1714,946)
    (1725,937)(1732,931)(1739,926)
    (1746,921)(1753,916)(1760,912)
    (1767,908)(1774,905)(1782,902)
    (1789,901)(1797,900)(1804,900)
    (1812,900)(1820,902)(1827,905)
    (1835,910)(1842,915)(1850,921)
    (1857,929)(1864,937)(1871,947)
    (1878,958)(1885,971)(1892,984)
    (1900,999)(1906,1015)(1913,1031)
    (1920,1050)(1927,1070)(1935,1093)
    (1943,1118)(1952,1146)(1961,1176)
    (1971,1210)(1981,1246)(1992,1286)
    (2004,1328)(2016,1373)(2029,1419)
    (2041,1467)(2054,1514)(2066,1559)
    (2077,1601)(2087,1639)(2095,1671)
    (2101,1696)(2106,1715)(2109,1727)
    (2111,1734)(2112,1737)
\path(2112,612)(2114,613)(2118,614)
    (2125,616)(2137,619)(2153,623)
    (2173,629)(2198,636)(2226,645)
    (2258,654)(2292,665)(2328,676)
    (2365,687)(2402,699)(2439,712)
    (2475,724)(2509,737)(2542,750)
    (2573,763)(2603,776)(2631,790)
    (2659,805)(2685,820)(2712,837)
    (2733,851)(2755,866)(2777,882)
    (2799,898)(2821,915)(2844,933)
    (2867,951)(2891,969)(2914,988)
    (2938,1008)(2962,1027)(2987,1047)
    (3011,1068)(3036,1088)(3060,1109)
    (3084,1130)(3109,1151)(3133,1172)
    (3157,1193)(3180,1214)(3203,1236)
    (3226,1257)(3248,1279)(3270,1301)
    (3291,1323)(3312,1345)(3331,1368)
    (3351,1390)(3369,1413)(3387,1437)
    (3404,1461)(3421,1486)(3437,1512)
    (3453,1539)(3470,1569)(3486,1599)
    (3502,1632)(3519,1667)(3536,1703)
    (3553,1741)(3570,1781)(3587,1822)
    (3605,1864)(3622,1906)(3640,1948)
    (3656,1990)(3672,2030)(3687,2069)
    (3701,2104)(3714,2137)(3725,2166)
    (3735,2192)(3743,2213)(3750,2230)
    (3755,2242)(3758,2252)(3760,2257)
    (3761,2261)(3762,2262)
\path(3762,2262)(3762,2261)(3764,2258)
    (3766,2252)(3769,2244)(3773,2232)
    (3780,2216)(3787,2196)(3797,2172)
    (3808,2144)(3821,2111)(3835,2075)
    (3851,2035)(3868,1993)(3886,1948)
    (3905,1902)(3925,1854)(3945,1806)
    (3966,1758)(3987,1710)(4008,1663)
    (4029,1618)(4050,1573)(4071,1531)
    (4092,1490)(4113,1452)(4133,1415)
    (4154,1381)(4175,1349)(4196,1318)
    (4218,1290)(4240,1262)(4263,1237)
    (4287,1212)(4312,1188)(4338,1165)
    (4366,1142)(4396,1119)(4427,1097)
    (4461,1074)(4498,1051)(4536,1028)
    (4577,1005)(4621,981)(4666,957)
    (4714,933)(4763,908)(4814,883)
    (4865,858)(4917,834)(4970,809)
    (5021,786)(5071,763)(5120,741)
    (5166,720)(5209,701)(5248,684)
    (5283,668)(5314,655)(5340,643)
    (5362,634)(5379,626)(5392,620)
    (5402,616)(5407,614)(5411,613)(5412,612)
\path(5412,2037)(5414,2035)(5419,2030)
    (5427,2022)(5439,2009)(5455,1992)
    (5475,1972)(5497,1948)(5521,1923)
    (5545,1897)(5569,1871)(5592,1846)
    (5614,1822)(5635,1799)(5653,1777)
    (5671,1757)(5686,1737)(5701,1719)
    (5714,1701)(5727,1684)(5738,1667)
    (5750,1649)(5763,1628)(5775,1607)
    (5788,1585)(5800,1563)(5812,1541)
    (5824,1519)(5836,1498)(5847,1477)
    (5859,1458)(5870,1440)(5881,1423)
    (5891,1409)(5901,1397)(5911,1388)
    (5921,1381)(5931,1376)(5940,1374)
    (5950,1374)(5959,1378)(5969,1383)
    (5980,1391)(5990,1400)(6002,1411)
    (6014,1423)(6026,1435)(6038,1447)
    (6051,1459)(6064,1469)(6077,1478)
    (6090,1485)(6103,1490)(6115,1491)
    (6127,1489)(6139,1484)(6151,1475)
    (6162,1462)(6169,1451)(6177,1438)
    (6184,1423)(6192,1405)(6200,1385)
    (6208,1362)(6216,1335)(6225,1305)
    (6235,1272)(6245,1235)(6255,1195)
    (6266,1151)(6277,1105)(6289,1055)
    (6301,1003)(6313,950)(6325,897)
    (6336,845)(6347,796)(6357,752)
    (6366,712)(6373,679)(6378,653)
    (6382,634)(6385,622)(6386,615)(6387,612)
\path(6387,612)(6387,615)(6388,621)
    (6390,633)(6392,651)(6395,676)
    (6399,708)(6405,746)(6411,791)
    (6417,841)(6425,894)(6432,950)
    (6440,1006)(6449,1063)(6457,1118)
    (6465,1171)(6473,1222)(6480,1270)
    (6488,1314)(6495,1355)(6502,1393)
    (6509,1427)(6516,1458)(6522,1486)
    (6529,1511)(6535,1533)(6542,1553)
    (6549,1571)(6555,1586)(6562,1599)
    (6575,1620)(6589,1636)(6604,1647)
    (6619,1654)(6635,1657)(6652,1658)
    (6669,1656)(6687,1653)(6705,1649)
    (6722,1646)(6739,1642)(6755,1641)
    (6770,1641)(6785,1645)(6799,1651)
    (6812,1662)(6820,1671)(6828,1681)
    (6835,1694)(6843,1708)(6850,1725)
    (6858,1743)(6865,1763)(6873,1785)
    (6880,1808)(6888,1834)(6895,1860)
    (6903,1889)(6910,1918)(6918,1949)
    (6925,1980)(6932,2012)(6939,2045)
    (6946,2078)(6953,2111)(6960,2143)
    (6967,2176)(6973,2209)(6980,2241)
    (6986,2273)(6993,2305)(7000,2337)
    (7005,2365)(7011,2392)(7017,2420)
    (7023,2448)(7030,2478)(7037,2508)
    (7044,2539)(7051,2571)(7059,2605)
    (7067,2641)(7076,2678)(7085,2717)
    (7095,2758)(7106,2801)(7116,2845)
    (7127,2890)(7138,2935)(7149,2979)
    (7159,3023)(7169,3064)(7179,3102)
    (7187,3136)(7194,3165)(7200,3190)
    (7205,3208)(7208,3222)(7210,3230)
    (7211,3235)(7212,3237)
\path(7212,612)(7213,614)(7214,618)
    (7217,626)(7221,639)(7227,657)
    (7235,681)(7245,711)(7257,747)
    (7271,789)(7287,836)(7304,888)
    (7323,944)(7342,1002)(7362,1062)
    (7382,1122)(7401,1180)(7420,1236)
    (7437,1288)(7453,1335)(7467,1377)
    (7479,1413)(7489,1443)(7497,1467)
    (7503,1485)(7507,1498)(7510,1506)
    (7511,1510)(7512,1512)
\path(7512,612)(7513,614)(7515,617)
    (7519,624)(7526,634)(7535,649)
    (7547,668)(7561,692)(7579,720)
    (7599,752)(7621,788)(7645,827)
    (7670,868)(7697,911)(7724,954)
    (7751,998)(7778,1041)(7804,1084)
    (7830,1126)(7855,1166)(7880,1205)
    (7903,1243)(7926,1279)(7948,1313)
    (7968,1347)(7989,1379)(8008,1410)
    (8027,1440)(8046,1470)(8064,1499)
    (8083,1527)(8101,1556)(8119,1584)
    (8137,1612)(8156,1642)(8176,1673)
    (8196,1703)(8216,1734)(8237,1765)
    (8258,1796)(8279,1828)(8301,1859)
    (8324,1891)(8347,1923)(8370,1954)
    (8393,1985)(8417,2015)(8442,2045)
    (8466,2074)(8490,2103)(8515,2130)
    (8540,2156)(8564,2181)(8589,2204)
    (8613,2226)(8637,2246)(8661,2265)
    (8684,2281)(8707,2296)(8730,2309)
    (8753,2321)(8775,2330)(8797,2338)
    (8819,2344)(8840,2347)(8862,2349)
    (8885,2350)(8908,2348)(8931,2343)
    (8955,2337)(8979,2329)(9003,2318)
    (9027,2306)(9052,2291)(9077,2275)
    (9103,2256)(9129,2235)(9155,2213)
    (9181,2189)(9207,2164)(9234,2137)
    (9260,2109)(9286,2079)(9311,2049)
    (9337,2018)(9362,1986)(9387,1954)
    (9411,1922)(9434,1889)(9457,1856)
    (9480,1824)(9502,1791)(9524,1759)
    (9545,1726)(9566,1694)(9587,1662)
    (9605,1634)(9623,1605)(9641,1577)
    (9660,1548)(9678,1518)(9697,1488)
    (9716,1457)(9735,1425)(9756,1391)
    (9776,1356)(9798,1320)(9821,1282)
    (9844,1243)(9869,1202)(9894,1159)
    (9920,1115)(9946,1070)(9973,1024)
    (10000,977)(10027,931)(10054,885)
    (10079,842)(10103,800)(10125,762)
    (10145,727)(10163,697)(10177,672)
    (10189,651)(10198,636)(10205,625)
    (10209,618)(10211,614)(10212,612)
\path(10212,3012)(10212,3013)(10214,3017)
    (10217,3028)(10224,3047)(10233,3074)
    (10246,3109)(10260,3149)(10276,3192)
    (10293,3235)(10310,3276)(10326,3314)
    (10342,3346)(10356,3373)(10370,3394)
    (10383,3410)(10394,3419)(10406,3422)
    (10416,3420)(10427,3412)(10437,3399)
    (10443,3389)(10450,3376)(10456,3361)
    (10463,3344)(10470,3324)(10477,3301)
    (10484,3275)(10491,3246)(10499,3214)
    (10507,3179)(10516,3140)(10525,3097)
    (10534,3052)(10544,3003)(10554,2950)
    (10564,2896)(10574,2839)(10585,2781)
    (10596,2722)(10606,2664)(10616,2608)
    (10625,2555)(10633,2507)(10641,2463)
    (10647,2427)(10652,2396)(10656,2373)
    (10659,2356)(10661,2345)(10662,2340)(10662,2337)
\path(10662,2337)(10662,2339)(10663,2345)
    (10664,2355)(10665,2370)(10667,2392)
    (10670,2420)(10673,2456)(10677,2498)
    (10681,2547)(10687,2602)(10692,2661)
    (10698,2723)(10705,2789)(10711,2855)
    (10718,2923)(10725,2989)(10732,3055)
    (10738,3118)(10745,3179)(10752,3237)
    (10758,3292)(10764,3344)(10770,3393)
    (10776,3438)(10782,3480)(10788,3519)
    (10793,3554)(10799,3587)(10804,3617)
    (10809,3645)(10815,3670)(10820,3693)
    (10826,3714)(10831,3732)(10837,3749)
    (10845,3770)(10853,3788)(10862,3802)
    (10871,3814)(10881,3821)(10891,3826)
    (10902,3828)(10914,3826)(10926,3821)
    (10940,3813)(10954,3803)(10970,3790)
    (10986,3774)(11002,3756)(11019,3737)
    (11035,3717)(11051,3697)(11065,3678)
    (11078,3660)(11089,3645)(11098,3633)
    (11104,3623)(11109,3617)(11111,3614)(11112,3612)
\put(792,3477){\makebox(0,0)[lb]{\smash{{{\SetFigFont{6}{7.2}{rm}$\varphi(x)$}}}}}
\end{picture}
}

%% file: fig3.eepic.tex
\setlength{\unitlength}{0.00041667in}
\begingroup\makeatletter\ifx\SetFigFont\undefined
\def\x#1#2#3#4#5#6#7\relax{\def\x{#1#2#3#4#5#6}}%
\expandafter\x\fmtname xxxxxx\relax \def\y{splain}%
\ifx\x\y   
\gdef\SetFigFont#1#2#3{%
  \ifnum #1<17\tiny\else \ifnum #1<20\small\else
  \ifnum #1<24\normalsize\else \ifnum #1<29\large\else
  \ifnum #1<34\Large\else \ifnum #1<41\LARGE\else
     \huge\fi\fi\fi\fi\fi\fi
  \csname #3\endcsname}%
\else
\gdef\SetFigFont#1#2#3{\begingroup
  \count@#1\relax \ifnum 25<\count@\count@25\fi
  \def\x{\endgroup\@setsize\SetFigFont{#2pt}}%
  \expandafter\x
    \csname \romannumeral\the\count@ pt\expandafter\endcsname
    \csname @\romannumeral\the\count@ pt\endcsname
  \csname #3\endcsname}%
\fi
\fi\endgroup
{\renewcommand{\dashlinestretch}{30}
\begin{picture}(12324,7440)(0,-10)
\put(8562,39){\makebox(0,0)[lb]{\smash{{{\SetFigFont{6}{7.2}{rm}$Q_{j+1}$}}}}}
\put(1587,39){\makebox(0,0)[lb]{\smash{{{\SetFigFont{6}{7.2}{rm}$Q_j$}}}}}
\put(2337,2499){\makebox(0,0)[lb]{\smash{{{\SetFigFont{6}{7.2}{rm}$x_{i-1}$}}}}}
\put(3537,2514){\makebox(0,0)[lb]{\smash{{{\SetFigFont{6}{7.2}{rm}$y_i^-$}}}}}
\put(6012,2499){\makebox(0,0)[lb]{\smash{{{\SetFigFont{6}{7.2}{rm}$x_i$}}}}}
\put(10812,2499){\makebox(0,0)[lb]{\smash{{{\SetFigFont{6}{7.2}{rm}$x_{i+1}$}}}}}
\put(1062,2499){\makebox(0,0)[lb]{\smash{{{\SetFigFont{6}{7.2}{rm}$x_{i-2}$}}}}}
\put(1737,2514){\makebox(0,0)[lb]{\smash{{{\SetFigFont{6}{7.2}{rm}$y_{i-1}$}}}}}
\put(9012,2514){\makebox(0,0)[lb]{\smash{{{\SetFigFont{6}{7.2}{rm}$y_{i+1}$ }}}}}
\put(6012,6444){\blacken\ellipse{90}{90}}
\put(6012,6444){\ellipse{90}{90}}
\put(10812,5244){\blacken\ellipse{90}{90}}
\put(10812,5244){\ellipse{90}{90}}
\put(1512,5994){\blacken\ellipse{90}{90}}
\put(1512,5994){\ellipse{90}{90}}
\put(2412,7044){\blacken\ellipse{90}{90}}
\put(2412,7044){\ellipse{90}{90}}
\path(612,2844)(11712,2844)
\path(612,444)(4212,444)
\path(612,444)(4212,444)
\blacken\path(4092.000,414.000)(4212.000,444.000)(4092.000,474.000)(4092.000,414.000)
\path(11712,444)(4212,444)
\path(11712,444)(4212,444)
\blacken\path(4332.000,474.000)(4212.000,444.000)(4332.000,414.000)(4332.000,474.000)
\path(3612,2844)(2412,7044)(1812,4944)
	(1512,5994)(912,3894)
\path(3613,4028)(6013,6458)
\path(6012,6444)(9012,3444)(10812,5244)(11712,4344)
\dashline{60.000}(9012,3444)(9012,2844)
\dashline{60.000}(10812,5169)(10812,2844)
\dashline{60.000}(6012,6444)(6012,2844)
\dashline{60.000}(3612,4044)(3612,2844)
\dashline{60.000}(2412,6969)(2412,2844)
\dashline{60.000}(1812,4944)(1812,2844)
\dashline{60.000}(1512,5919)(1512,2844)
\path(1962,2844)(1962,4944)
\path(1962,2844)(1962,4944)
\blacken\path(1992.000,4824.000)(1962.000,4944.000)(1932.000,4824.000)(1992.000,4824.000)
\path(1962,4944)(1962,2844)
\path(1962,4944)(1962,2844)
\blacken\path(1932.000,2964.000)(1962.000,2844.000)(1992.000,2964.000)(1932.000,2964.000)
\dottedline{45}(1812,4944)(1962,4944)
\path(8862,3444)(8862,2844)
\path(8862,3444)(8862,2844)
\blacken\path(8832.000,2964.000)(8862.000,2844.000)(8892.000,2964.000)(8832.000,2964.000)
\path(8862,2844)(8862,3444)
\path(8862,2844)(8862,3444)
\blacken\path(8892.000,3324.000)(8862.000,3444.000)(8832.000,3324.000)(8892.000,3324.000)
\dottedline{45}(9012,3444)(8862,3444)
\dashline{60.000}(4212,6144)(4212,444)
\dottedline{45}(3612,4044)(3162,4044)
\path(3162,2844)(3162,4044)
\path(3162,2844)(3162,4044)
\blacken\path(3192.000,3924.000)(3162.000,4044.000)(3132.000,3924.000)(3192.000,3924.000)
\path(3162,4044)(3162,2844)
\path(3162,4044)(3162,2844)
\blacken\path(3132.000,2964.000)(3162.000,2844.000)(3192.000,2964.000)(3132.000,2964.000)
\dottedline{90}(12,444)(612,444)
\dottedline{90}(11712,444)(12312,444)
\thicklines
\path(612,2244)(613,2244)(616,2243)
	(621,2241)(630,2238)(643,2234)
	(659,2229)(681,2222)(708,2214)
	(739,2204)(777,2192)(819,2178)
	(867,2163)(920,2147)(978,2129)
	(1040,2110)(1105,2090)(1174,2068)
	(1246,2046)(1320,2024)(1395,2001)
	(1472,1978)(1549,1955)(1626,1932)
	(1704,1909)(1780,1886)(1856,1864)
	(1930,1843)(2003,1822)(2074,1802)
	(2144,1783)(2211,1764)(2277,1746)
	(2342,1729)(2404,1713)(2464,1698)
	(2523,1683)(2580,1669)(2635,1656)
	(2689,1644)(2742,1632)(2793,1621)
	(2843,1611)(2892,1601)(2940,1592)
	(2988,1584)(3034,1576)(3081,1569)
	(3126,1562)(3172,1555)(3217,1549)
	(3262,1544)(3311,1539)(3359,1533)
	(3408,1529)(3457,1525)(3506,1521)
	(3555,1518)(3605,1516)(3656,1514)
	(3708,1512)(3760,1511)(3814,1511)
	(3869,1511)(3925,1511)(3983,1512)
	(4043,1514)(4104,1516)(4167,1518)
	(4232,1521)(4299,1525)(4368,1528)
	(4438,1533)(4510,1537)(4583,1542)
	(4657,1548)(4733,1553)(4809,1559)
	(4885,1566)(4960,1572)(5035,1579)
	(5108,1585)(5180,1591)(5248,1598)
	(5313,1604)(5374,1610)(5431,1615)
	(5483,1620)(5529,1625)(5570,1629)
	(5605,1633)(5635,1636)(5659,1638)
	(5678,1640)(5692,1642)(5701,1643)
	(5707,1643)(5711,1644)(5712,1644)
\path(5712,1644)(5714,1645)(5719,1649)
	(5727,1654)(5740,1663)(5758,1676)
	(5782,1692)(5811,1712)(5846,1734)
	(5885,1759)(5928,1787)(5973,1815)
	(6021,1844)(6070,1873)(6119,1901)
	(6167,1928)(6215,1953)(6261,1976)
	(6306,1997)(6349,2016)(6390,2031)
	(6429,2045)(6466,2055)(6502,2063)
	(6536,2068)(6569,2071)(6601,2071)
	(6632,2068)(6662,2063)(6691,2055)
	(6720,2044)(6750,2031)(6773,2020)
	(6796,2006)(6819,1992)(6843,1975)
	(6867,1958)(6892,1939)(6917,1919)
	(6942,1898)(6969,1876)(6995,1852)
	(7023,1828)(7052,1802)(7081,1776)
	(7111,1749)(7142,1722)(7174,1694)
	(7207,1665)(7240,1637)(7275,1608)
	(7310,1579)(7346,1551)(7384,1523)
	(7421,1495)(7460,1467)(7500,1441)
	(7540,1415)(7581,1390)(7623,1366)
	(7665,1343)(7709,1321)(7753,1300)
	(7797,1281)(7843,1263)(7889,1247)
	(7936,1232)(7984,1219)(8033,1207)
	(8083,1197)(8134,1188)(8187,1181)
	(8225,1178)(8265,1175)(8305,1173)
	(8346,1172)(8389,1172)(8432,1173)
	(8477,1175)(8524,1178)(8572,1182)
	(8621,1187)(8673,1193)(8726,1201)
	(8781,1209)(8838,1219)(8898,1230)
	(8960,1242)(9024,1255)(9090,1270)
	(9159,1285)(9231,1302)(9305,1320)
	(9381,1339)(9460,1360)(9542,1381)
	(9625,1404)(9711,1427)(9799,1451)
	(9888,1477)(9979,1503)(10070,1529)
	(10163,1556)(10256,1584)(10348,1611)
	(10440,1639)(10530,1666)(10618,1693)
	(10704,1720)(10786,1745)(10865,1770)
	(10940,1793)(11009,1815)(11074,1836)
	(11133,1854)(11186,1871)(11233,1886)
	(11274,1899)(11309,1910)(11338,1920)
	(11361,1927)(11379,1933)(11393,1938)
	(11402,1941)(11408,1943)(11411,1944)(11412,1944)
\path(11412,1944)(11415,1945)(11421,1949)
	(11432,1954)(11449,1962)(11471,1973)
	(11498,1987)(11529,2003)(11562,2019)
	(11595,2035)(11626,2051)(11653,2065)
	(11675,2076)(11692,2084)(11703,2089)
	(11709,2093)(11712,2094)
\path(3312,5244)(3312,5243)(3313,5240)
	(3315,5235)(3318,5227)(3323,5215)
	(3329,5200)(3336,5179)(3346,5154)
	(3357,5123)(3371,5087)(3387,5045)
	(3404,4998)(3425,4944)(3447,4884)
	(3471,4819)(3498,4749)(3526,4673)
	(3556,4593)(3588,4508)(3621,4420)
	(3655,4328)(3690,4235)(3726,4140)
	(3762,4044)(3798,3948)(3834,3853)
	(3869,3760)(3903,3668)(3936,3580)
	(3968,3495)(3998,3415)(4026,3339)
	(4053,3269)(4077,3204)(4099,3144)
	(4120,3090)(4137,3043)(4153,3001)
	(4167,2965)(4178,2934)(4188,2909)
	(4195,2888)(4201,2873)(4206,2861)
	(4209,2853)(4211,2848)(4212,2845)(4212,2844)
\path(4212,2844)(4213,2842)(4216,2838)
	(4221,2831)(4228,2820)(4239,2804)
	(4253,2784)(4270,2758)(4291,2729)
	(4314,2695)(4339,2658)(4366,2619)
	(4394,2578)(4423,2536)(4453,2494)
	(4482,2452)(4511,2411)(4539,2371)
	(4567,2332)(4593,2295)(4619,2259)
	(4644,2225)(4667,2193)(4690,2162)
	(4713,2132)(4735,2103)(4756,2075)
	(4777,2048)(4798,2022)(4819,1996)
	(4841,1970)(4862,1944)(4881,1921)
	(4901,1898)(4921,1875)(4941,1852)
	(4962,1829)(4984,1805)(5007,1782)
	(5030,1758)(5054,1734)(5079,1711)
	(5106,1687)(5133,1664)(5161,1640)
	(5190,1617)(5220,1594)(5252,1572)
	(5284,1550)(5317,1529)(5352,1508)
	(5387,1488)(5423,1468)(5460,1450)
	(5498,1432)(5537,1415)(5576,1399)
	(5616,1385)(5657,1371)(5699,1358)
	(5741,1346)(5785,1336)(5828,1326)
	(5873,1318)(5919,1310)(5965,1304)
	(6013,1298)(6062,1294)(6098,1292)
	(6135,1290)(6173,1288)(6212,1287)
	(6251,1287)(6292,1287)(6334,1287)
	(6377,1289)(6420,1290)(6465,1292)
	(6512,1295)(6559,1298)(6607,1301)
	(6656,1305)(6707,1310)(6758,1315)
	(6811,1320)(6864,1326)(6919,1333)
	(6974,1339)(7030,1346)(7087,1354)
	(7144,1362)(7202,1370)(7261,1379)
	(7319,1388)(7379,1397)(7438,1406)
	(7498,1416)(7558,1426)(7618,1436)
	(7678,1446)(7738,1457)(7798,1468)
	(7858,1478)(7917,1489)(7976,1500)
	(8035,1511)(8093,1522)(8151,1533)
	(8209,1544)(8266,1556)(8323,1567)
	(8379,1578)(8435,1589)(8491,1600)
	(8547,1611)(8602,1622)(8657,1633)
	(8712,1644)(8763,1654)(8814,1664)
	(8865,1675)(8916,1685)(8968,1695)
	(9020,1706)(9073,1716)(9126,1727)
	(9181,1738)(9236,1749)(9292,1760)
	(9349,1771)(9408,1783)(9468,1795)
	(9530,1808)(9593,1820)(9658,1833)
	(9724,1846)(9793,1860)(9863,1874)
	(9936,1889)(10010,1904)(10085,1919)
	(10163,1934)(10242,1950)(10323,1966)
	(10405,1983)(10487,1999)(10571,2016)
	(10654,2032)(10738,2049)(10820,2066)
	(10902,2082)(10982,2098)(11060,2114)
	(11136,2129)(11208,2143)(11276,2157)
	(11340,2170)(11399,2181)(11453,2192)
	(11502,2202)(11546,2211)(11584,2218)
	(11616,2225)(11643,2230)(11665,2235)
	(11681,2238)(11694,2240)(11703,2242)
	(11708,2243)(11711,2244)(11712,2244)
\path(1812,6069)(1813,6072)(1815,6077)
	(1819,6088)(1825,6104)(1833,6125)
	(1843,6151)(1855,6181)(1868,6215)
	(1882,6250)(1896,6286)(1910,6322)
	(1924,6356)(1937,6389)(1949,6420)
	(1961,6449)(1972,6476)(1983,6500)
	(1992,6523)(2002,6545)(2011,6565)
	(2020,6583)(2028,6602)(2037,6619)
	(2048,6642)(2060,6664)(2072,6685)
	(2084,6706)(2096,6727)(2109,6748)
	(2122,6768)(2135,6788)(2148,6807)
	(2161,6825)(2173,6842)(2185,6858)
	(2197,6873)(2209,6887)(2219,6899)
	(2230,6911)(2240,6922)(2250,6931)
	(2262,6943)(2274,6954)(2286,6964)
	(2299,6974)(2313,6984)(2328,6995)
	(2344,7005)(2361,7015)(2377,7025)
	(2392,7033)(2402,7039)(2409,7042)(2412,7044)
\path(912,5544)(914,5546)(919,5551)
	(928,5560)(941,5572)(957,5588)
	(976,5607)(998,5627)(1020,5649)
	(1043,5670)(1065,5691)(1085,5710)
	(1105,5728)(1123,5745)(1139,5760)
	(1155,5773)(1170,5786)(1184,5798)
	(1198,5808)(1212,5819)(1227,5830)
	(1243,5841)(1259,5852)(1275,5863)
	(1292,5873)(1309,5884)(1326,5894)
	(1344,5904)(1362,5914)(1380,5923)
	(1398,5932)(1415,5941)(1432,5949)
	(1449,5956)(1465,5963)(1481,5969)
	(1497,5976)(1512,5981)(1529,5988)
	(1547,5994)(1565,6000)(1585,6006)
	(1606,6013)(1629,6019)(1654,6027)
	(1681,6034)(1709,6042)(1737,6049)
	(1762,6056)(1783,6062)(1799,6066)
	(1808,6068)(1811,6069)(1812,6069)
\path(2412,7044)(2414,7043)(2419,7040)
	(2428,7036)(2441,7029)(2457,7020)
	(2476,7008)(2496,6995)(2517,6981)
	(2537,6966)(2556,6950)(2574,6933)
	(2592,6914)(2609,6895)(2625,6873)
	(2641,6849)(2658,6823)(2675,6794)
	(2685,6775)(2696,6755)(2707,6734)
	(2718,6712)(2730,6688)(2742,6663)
	(2754,6638)(2767,6610)(2780,6582)
	(2793,6553)(2806,6523)(2820,6492)
	(2834,6460)(2848,6427)(2862,6394)
	(2876,6361)(2890,6328)(2904,6294)
	(2918,6261)(2931,6228)(2944,6196)
	(2957,6164)(2970,6132)(2982,6101)
	(2994,6071)(3006,6042)(3017,6013)
	(3028,5985)(3039,5958)(3050,5931)
	(3061,5901)(3073,5872)(3084,5843)
	(3095,5814)(3107,5785)(3118,5755)
	(3130,5725)(3142,5694)(3154,5662)
	(3167,5628)(3180,5593)(3194,5557)
	(3207,5521)(3222,5483)(3236,5446)
	(3249,5410)(3262,5375)(3274,5344)
	(3285,5316)(3294,5292)(3301,5273)
	(3306,5260)(3309,5251)(3311,5246)(3312,5244)
\path(6012,6444)(6014,6443)(6017,6441)
	(6024,6438)(6034,6432)(6048,6424)
	(6068,6413)(6092,6400)(6120,6384)
	(6154,6365)(6191,6344)(6232,6322)
	(6275,6297)(6321,6271)(6368,6244)
	(6417,6217)(6465,6189)(6513,6161)
	(6560,6133)(6607,6106)(6652,6079)
	(6696,6053)(6738,6027)(6779,6002)
	(6819,5978)(6857,5954)(6894,5930)
	(6930,5907)(6965,5884)(7000,5861)
	(7034,5838)(7067,5815)(7100,5791)
	(7133,5768)(7166,5744)(7200,5719)
	(7230,5696)(7260,5673)(7291,5650)
	(7322,5626)(7354,5601)(7386,5576)
	(7418,5550)(7451,5524)(7485,5497)
	(7519,5469)(7553,5441)(7588,5412)
	(7623,5383)(7658,5354)(7694,5324)
	(7730,5294)(7766,5263)(7801,5232)
	(7837,5202)(7873,5171)(7909,5140)
	(7944,5109)(7979,5079)(8013,5048)
	(8047,5018)(8081,4989)(8114,4959)
	(8146,4931)(8178,4902)(8208,4875)
	(8238,4848)(8268,4821)(8296,4795)
	(8324,4770)(8351,4746)(8377,4722)
	(8402,4698)(8427,4675)(8451,4653)
	(8475,4631)(8505,4603)(8535,4576)
	(8563,4549)(8591,4523)(8619,4497)
	(8646,4471)(8674,4445)(8701,4419)
	(8729,4392)(8757,4365)(8785,4338)
	(8815,4310)(8844,4281)(8875,4251)
	(8905,4221)(8936,4191)(8966,4162)
	(8996,4132)(9025,4104)(9052,4078)
	(9076,4053)(9098,4032)(9117,4013)
	(9132,3998)(9144,3987)(9152,3979)
	(9158,3973)(9161,3970)(9162,3969)
\path(10812,5244)(10815,5244)(10823,5246)
	(10836,5247)(10855,5250)(10879,5253)
	(10909,5257)(10941,5261)(10974,5265)
	(11008,5268)(11041,5271)(11072,5274)
	(11101,5276)(11128,5277)(11153,5277)
	(11177,5277)(11199,5276)(11220,5274)
	(11241,5272)(11262,5269)(11283,5265)
	(11304,5261)(11325,5255)(11347,5249)
	(11371,5241)(11396,5232)(11423,5222)
	(11452,5210)(11483,5197)(11516,5183)
	(11550,5169)(11583,5154)(11615,5139)
	(11645,5126)(11669,5114)(11688,5105)
	(11701,5099)(11709,5096)(11712,5094)
\path(4212,6144)(4214,6146)(4220,6149)
	(4229,6155)(4244,6165)(4264,6178)
	(4289,6194)(4318,6213)(4352,6234)
	(4388,6257)(4426,6281)(4464,6305)
	(4502,6329)(4539,6352)(4574,6374)
	(4608,6395)(4640,6414)(4670,6432)
	(4698,6448)(4725,6463)(4750,6477)
	(4774,6489)(4797,6501)(4819,6512)
	(4841,6522)(4862,6531)(4887,6542)
	(4912,6552)(4937,6562)(4962,6570)
	(4988,6579)(5014,6587)(5040,6594)
	(5067,6600)(5094,6606)(5121,6611)
	(5148,6616)(5175,6620)(5202,6623)
	(5228,6625)(5254,6627)(5280,6627)
	(5305,6628)(5330,6627)(5354,6626)
	(5378,6624)(5401,6622)(5425,6619)
	(5446,6616)(5468,6612)(5490,6608)
	(5513,6603)(5538,6597)(5563,6591)
	(5590,6584)(5619,6575)(5650,6566)
	(5682,6556)(5717,6545)(5753,6533)
	(5790,6521)(5827,6509)(5863,6496)
	(5897,6484)(5929,6474)(5955,6464)
	(5977,6457)(5993,6451)(6003,6447)
	(6009,6445)(6012,6444)
\path(9162,3969)(9164,3971)(9168,3975)
	(9175,3982)(9186,3992)(9201,4008)
	(9221,4028)(9246,4052)(9275,4081)
	(9307,4113)(9343,4148)(9381,4186)
	(9420,4225)(9461,4265)(9501,4305)
	(9541,4344)(9580,4382)(9618,4419)
	(9654,4455)(9689,4488)(9721,4520)
	(9752,4550)(9781,4578)(9809,4604)
	(9835,4629)(9860,4652)(9883,4674)
	(9905,4695)(9927,4715)(9947,4733)
	(9967,4752)(9987,4769)(10012,4791)
	(10037,4813)(10062,4834)(10086,4855)
	(10111,4875)(10135,4894)(10158,4913)
	(10182,4932)(10205,4950)(10228,4967)
	(10251,4984)(10273,5000)(10295,5015)
	(10316,5030)(10336,5043)(10356,5056)
	(10375,5068)(10393,5079)(10410,5089)
	(10426,5099)(10442,5108)(10458,5116)
	(10473,5124)(10487,5131)(10506,5141)
	(10525,5149)(10543,5158)(10563,5166)
	(10583,5173)(10604,5181)(10627,5189)
	(10651,5197)(10677,5205)(10703,5213)
	(10729,5220)(10753,5227)(10774,5233)
	(10791,5238)(10802,5241)(10809,5243)(10812,5244)
\put(1287,1569){\makebox(0,0)[lb]{\smash{{{\SetFigFont{6}{7.2}{rm}$g_{j-1}(x)$}}}}}
\put(7587,5844){\makebox(0,0)[lb]{\smash{{{\SetFigFont{6}{7.2}{rm}$g_{j+1}(x)$}}}}}
\put(8262,1869){\makebox(0,0)[lb]{\smash{{{\SetFigFont{6}{7.2}{rm}$g_j(x)$}}}}}
\put(4287,2919){\makebox(0,0)[lb]{\smash{{{\SetFigFont{6}{7.2}{rm}$\omega$}}}}}
\put(987,6294){\makebox(0,0)[lb]{\smash{{{\SetFigFont{6}{7.2}{rm}$z_{i-2}$}}}}}
\put(2037,3969){\makebox(0,0)[lb]{\smash{{{\SetFigFont{6}{7.2}{rm}$R_{i-1}$}}}}}
\put(7962,3069){\makebox(0,0)[lb]{\smash{{{\SetFigFont{6}{7.2}{rm}$R_{i+1}$}}}}}
\put(4662,4794){\makebox(0,0)[lb]{\smash{{{\SetFigFont{6}{7.2}{rm}$\psi_i(x)$}}}}}
\put(10812,5469){\makebox(0,0)[lb]{\smash{{{\SetFigFont{6}{7.2}{rm}$z_{i+1}$}}}}}
\put(6012,6594){\makebox(0,0)[lb]{\smash{{{\SetFigFont{6}{7.2}{rm}$z_i$}}}}}
\put(2412,7269){\makebox(0,0)[lb]{\smash{{{\SetFigFont{6}{7.2}{rm}$z_{i-1}$}}}}}
\put(2712,3369){\makebox(0,0)[lb]{\smash{{{\SetFigFont{6}{7.2}{rm}$R_i$}}}}}
\put(6537,4719){\makebox(0,0)[lb]{\smash{{{\SetFigFont{6}{7.2}{rm}$\psi_{i+1}(x)$}}}}}
\end{picture}
}

%% file: fig4.eepic.tex
\setlength{\unitlength}{0.00050000in}
\begingroup\makeatletter\ifx\SetFigFont\undefined%
\gdef\SetFigFont#1#2#3#4#5{%
  \reset@font\fontsize{#1}{#2pt}%
  \fontfamily{#3}\fontseries{#4}\fontshape{#5}%
  \selectfont}%
\fi\endgroup%
{\renewcommand{\dashlinestretch}{30}
\begin{picture}(12324,6634)(0,-10)
\put(1737,58){\makebox(0,0)[lb]{\smash{{{\SetFigFont{6}{7.2}{rm}$Q_{j+2}$}}}}}
\put(3687,58){\makebox(0,0)[lb]{\smash{{{\SetFigFont{6}{7.2}{rm}$Q_{j+1}$}}}}}
\put(7737,58){\makebox(0,0)[lb]{\smash{{{\SetFigFont{6}{7.2}{rm}$Q_j$}}}}}
\thicklines
\path(612,1663)(614,1662)(618,1660)
	(626,1656)(638,1651)(655,1643)
	(677,1632)(705,1619)(738,1604)
	(776,1587)(818,1569)(863,1549)
	(911,1528)(961,1507)(1012,1485)
	(1064,1465)(1115,1445)(1166,1425)
	(1215,1408)(1263,1391)(1310,1376)
	(1355,1362)(1398,1350)(1440,1339)
	(1480,1330)(1520,1322)(1558,1316)
	(1595,1312)(1632,1309)(1668,1307)
	(1704,1306)(1740,1307)(1776,1310)
	(1812,1313)(1845,1317)(1878,1322)
	(1911,1328)(1946,1335)(1980,1342)
	(2016,1350)(2052,1359)(2089,1369)
	(2126,1379)(2165,1390)(2204,1401)
	(2244,1413)(2285,1425)(2326,1437)
	(2368,1450)(2410,1463)(2453,1475)
	(2496,1488)(2540,1501)(2584,1513)
	(2628,1526)(2672,1538)(2716,1549)
	(2760,1560)(2804,1570)(2848,1580)
	(2892,1589)(2935,1597)(2978,1604)
	(3020,1611)(3062,1616)(3104,1621)
	(3145,1624)(3187,1627)(3227,1628)
	(3268,1628)(3309,1628)(3350,1625)
	(3386,1623)(3423,1618)(3461,1613)
	(3499,1607)(3537,1599)(3576,1590)
	(3616,1580)(3657,1568)(3699,1554)
	(3743,1540)(3788,1523)(3834,1505)
	(3883,1486)(3933,1464)(3985,1441)
	(4038,1417)(4094,1390)(4152,1362)
	(4211,1333)(4272,1302)(4334,1270)
	(4397,1237)(4461,1203)(4525,1169)
	(4588,1135)(4650,1101)(4710,1067)
	(4768,1035)(4823,1004)(4874,975)
	(4920,949)(4962,925)(4998,904)
	(5029,886)(5055,872)(5075,860)
	(5090,851)(5100,845)(5107,841)
	(5110,839)(5112,838)
\path(5112,838)(5114,838)(5118,839)
	(5125,840)(5137,842)(5153,845)
	(5175,849)(5203,854)(5236,860)
	(5274,866)(5318,873)(5366,881)
	(5418,890)(5473,899)(5530,908)
	(5588,917)(5648,926)(5707,934)
	(5766,943)(5824,951)(5880,958)
	(5935,965)(5989,972)(6040,977)
	(6090,982)(6137,987)(6183,990)
	(6228,993)(6270,996)(6312,997)
	(6352,999)(6391,999)(6429,999)
	(6466,998)(6503,996)(6540,994)
	(6576,991)(6612,988)(6648,984)
	(6685,980)(6721,975)(6758,969)
	(6795,964)(6833,957)(6871,951)
	(6909,944)(6949,936)(6988,929)
	(7028,921)(7069,913)(7110,905)
	(7151,897)(7193,889)(7235,881)
	(7277,874)(7320,866)(7362,859)
	(7404,852)(7447,845)(7489,839)
	(7531,833)(7573,828)(7614,824)
	(7655,820)(7696,817)(7736,815)
	(7775,813)(7815,812)(7853,812)
	(7891,813)(7929,815)(7966,817)
	(8003,821)(8039,826)(8076,831)
	(8112,838)(8146,845)(8181,854)
	(8216,863)(8251,873)(8286,885)
	(8322,897)(8358,910)(8395,924)
	(8432,938)(8470,954)(8508,970)
	(8547,987)(8587,1004)(8627,1022)
	(8667,1041)(8708,1060)(8750,1079)
	(8791,1098)(8834,1118)(8876,1138)
	(8918,1158)(8961,1178)(9004,1197)
	(9047,1216)(9089,1235)(9132,1254)
	(9174,1272)(9217,1289)(9259,1306)
	(9301,1322)(9343,1338)(9384,1352)
	(9426,1366)(9467,1379)(9508,1391)
	(9548,1403)(9589,1413)(9630,1422)
	(9671,1431)(9712,1438)(9750,1444)
	(9788,1449)(9826,1453)(9865,1456)
	(9905,1458)(9946,1459)(9987,1460)
	(10030,1459)(10074,1458)(10120,1455)
	(10167,1451)(10216,1447)(10267,1441)
	(10320,1435)(10374,1427)(10431,1418)
	(10490,1409)(10551,1398)(10614,1387)
	(10679,1375)(10746,1362)(10814,1348)
	(10883,1333)(10953,1318)(11023,1303)
	(11093,1287)(11162,1272)(11229,1256)
	(11293,1241)(11355,1226)(11413,1212)
	(11466,1199)(11515,1188)(11558,1177)
	(11596,1168)(11627,1160)(11654,1153)
	(11674,1148)(11689,1144)(11700,1141)
	(11707,1139)(11710,1138)(11712,1138)
\thinlines
\put(9462,5488){\blacken\ellipse{90}{90}}
\put(9462,5488){\ellipse{90}{90}}
\put(1812,5908){\blacken\ellipse{90}{90}}
\put(1812,5908){\ellipse{90}{90}}
\path(612,2788)(11712,2788)
\path(912,4888)(1812,5938)(4512,2788)(5112,2038)
\path(5112,2038)(5412,1663)
\dashline{60.000}(1812,2788)(1812,5938)
\dashline{60.000}(5412,2788)(5412,1663)
\dashline{60.000}(3612,6088)(3612,388)
\dashline{60.000}(4212,4288)(4212,388)
\path(11712,388)(4212,388)
\path(11712,388)(4212,388)
\blacken\path(4332.000,418.000)(4212.000,388.000)(4332.000,358.000)(4332.000,418.000)
\path(3612,388)(4212,388)
\path(3612,388)(4212,388)
\blacken\path(4092.000,358.000)(4212.000,388.000)(4092.000,418.000)(4092.000,358.000)
\path(4212,388)(3612,388)
\path(4212,388)(3612,388)
\blacken\path(3732.000,418.000)(3612.000,388.000)(3732.000,358.000)(3732.000,418.000)
\path(612,388)(3612,388)
\path(612,388)(3612,388)
\blacken\path(3492.000,358.000)(3612.000,388.000)(3492.000,418.000)(3492.000,358.000)
\path(5862,2788)(5862,1663)
\path(5862,2788)(5862,1663)
\blacken\path(5832.000,1783.000)(5862.000,1663.000)(5892.000,1783.000)(5832.000,1783.000)
\path(5862,1663)(5862,2788)
\path(5862,1663)(5862,2788)
\blacken\path(5892.000,2668.000)(5862.000,2788.000)(5832.000,2668.000)(5892.000,2668.000)
\dottedline{45}(5412,1663)(6012,1663)
\path(5412,2788)(9462,5488)(11637,4288)
\dashline{60.000}(9462,5488)(9462,2788)
\dottedline{90}(12,388)(612,388)
\dottedline{90}(11712,388)(12312,388)
\thicklines
\path(1812,5938)(1815,5939)(1821,5942)
	(1831,5947)(1847,5955)(1869,5966)
	(1896,5979)(1927,5994)(1962,6011)
	(1999,6029)(2036,6047)(2073,6064)
	(2110,6082)(2145,6098)(2178,6114)
	(2210,6128)(2239,6141)(2267,6154)
	(2293,6165)(2318,6176)(2342,6186)
	(2366,6195)(2389,6204)(2412,6213)
	(2435,6221)(2458,6230)(2482,6238)
	(2505,6246)(2530,6254)(2555,6261)
	(2580,6268)(2606,6275)(2632,6282)
	(2658,6288)(2685,6294)(2712,6299)
	(2739,6304)(2766,6308)(2792,6311)
	(2818,6314)(2844,6317)(2870,6318)
	(2894,6319)(2919,6319)(2942,6319)
	(2966,6317)(2989,6316)(3012,6313)
	(3035,6310)(3058,6306)(3082,6301)
	(3106,6295)(3131,6288)(3157,6280)
	(3185,6271)(3214,6260)(3246,6248)
	(3279,6235)(3314,6221)(3351,6206)
	(3388,6190)(3425,6173)(3462,6157)
	(3497,6141)(3528,6127)(3555,6115)
	(3577,6105)(3593,6097)(3603,6092)
	(3609,6089)(3612,6088)
\path(612,2188)(614,2187)(618,2185)
	(625,2180)(636,2173)(652,2164)
	(673,2151)(699,2136)(729,2118)
	(763,2098)(801,2077)(841,2054)
	(883,2031)(926,2007)(970,1984)
	(1013,1961)(1055,1939)(1097,1919)
	(1137,1900)(1176,1882)(1214,1866)
	(1250,1852)(1284,1839)(1318,1827)
	(1350,1818)(1382,1809)(1412,1803)
	(1443,1797)(1472,1793)(1502,1790)
	(1532,1789)(1562,1788)(1590,1788)
	(1619,1790)(1648,1792)(1678,1795)
	(1708,1799)(1739,1804)(1771,1810)
	(1803,1817)(1835,1825)(1868,1834)
	(1902,1843)(1936,1853)(1971,1864)
	(2006,1876)(2041,1889)(2076,1902)
	(2112,1916)(2148,1930)(2183,1945)
	(2218,1960)(2253,1976)(2288,1992)
	(2322,2008)(2356,2024)(2389,2041)
	(2421,2057)(2453,2074)(2485,2090)
	(2516,2107)(2546,2123)(2576,2139)
	(2605,2156)(2634,2172)(2662,2188)
	(2692,2205)(2722,2222)(2752,2240)
	(2781,2258)(2812,2276)(2842,2294)
	(2874,2313)(2906,2333)(2940,2354)
	(2974,2375)(3010,2398)(3048,2422)
	(3087,2447)(3127,2472)(3169,2499)
	(3211,2527)(3254,2555)(3298,2583)
	(3341,2611)(3383,2638)(3423,2664)
	(3461,2689)(3495,2711)(3525,2731)
	(3551,2748)(3572,2761)(3588,2772)
	(3599,2779)(3606,2784)(3610,2787)(3612,2788)
\path(3612,2788)(3613,2790)(3614,2795)
	(3616,2805)(3620,2819)(3625,2840)
	(3632,2866)(3641,2898)(3650,2936)
	(3662,2979)(3674,3026)(3687,3076)
	(3701,3128)(3715,3180)(3729,3233)
	(3743,3285)(3756,3336)(3770,3384)
	(3783,3431)(3795,3475)(3807,3517)
	(3818,3556)(3829,3593)(3839,3627)
	(3849,3659)(3859,3689)(3868,3718)
	(3877,3744)(3886,3769)(3895,3793)
	(3903,3816)(3912,3838)(3925,3868)
	(3937,3897)(3950,3925)(3963,3951)
	(3977,3977)(3992,4002)(4008,4028)
	(4025,4053)(4042,4079)(4061,4105)
	(4081,4131)(4102,4157)(4122,4183)
	(4142,4207)(4160,4229)(4176,4248)
	(4190,4263)(4200,4274)(4207,4282)
	(4210,4286)(4212,4288)
\path(612,2488)(614,2487)(618,2485)
	(625,2481)(636,2474)(652,2465)
	(673,2453)(699,2439)(729,2422)
	(763,2402)(801,2381)(841,2359)
	(883,2335)(926,2311)(970,2287)
	(1013,2264)(1055,2240)(1097,2218)
	(1137,2196)(1176,2175)(1214,2155)
	(1250,2137)(1284,2119)(1318,2102)
	(1350,2086)(1382,2070)(1412,2056)
	(1443,2041)(1472,2028)(1502,2014)
	(1532,2001)(1562,1988)(1592,1975)
	(1623,1962)(1654,1950)(1686,1937)
	(1718,1925)(1751,1912)(1784,1900)
	(1819,1887)(1854,1875)(1889,1863)
	(1925,1851)(1962,1839)(1999,1827)
	(2037,1816)(2074,1805)(2112,1795)
	(2150,1784)(2187,1775)(2225,1766)
	(2262,1757)(2299,1749)(2335,1741)
	(2370,1735)(2405,1728)(2440,1723)
	(2473,1718)(2506,1713)(2538,1710)
	(2570,1706)(2601,1704)(2632,1702)
	(2662,1700)(2694,1700)(2726,1700)
	(2758,1700)(2789,1701)(2821,1703)
	(2854,1706)(2886,1710)(2919,1714)
	(2951,1719)(2984,1725)(3017,1732)
	(3050,1740)(3083,1748)(3116,1757)
	(3148,1767)(3180,1778)(3212,1790)
	(3243,1802)(3274,1815)(3303,1828)
	(3333,1842)(3361,1856)(3389,1871)
	(3415,1887)(3441,1903)(3467,1919)
	(3491,1935)(3515,1952)(3539,1970)
	(3562,1988)(3583,2005)(3605,2024)
	(3626,2042)(3647,2062)(3668,2083)
	(3689,2105)(3711,2128)(3734,2153)
	(3757,2179)(3780,2207)(3805,2237)
	(3831,2269)(3857,2302)(3885,2338)
	(3913,2375)(3942,2413)(3971,2452)
	(4001,2492)(4030,2532)(4058,2571)
	(4085,2609)(4110,2644)(4133,2676)
	(4154,2705)(4171,2729)(4185,2749)
	(4196,2765)(4203,2775)(4208,2782)
	(4211,2786)(4212,2788)
\path(612,4888)(614,4890)(617,4893)
	(623,4900)(633,4911)(646,4926)
	(663,4944)(683,4967)(706,4993)
	(732,5021)(759,5051)(787,5082)
	(816,5114)(845,5145)(873,5176)
	(901,5205)(927,5233)(952,5260)
	(977,5285)(1000,5309)(1022,5331)
	(1043,5352)(1064,5372)(1084,5392)
	(1103,5410)(1123,5428)(1142,5446)
	(1162,5463)(1182,5480)(1202,5497)
	(1222,5514)(1243,5531)(1265,5549)
	(1288,5567)(1312,5585)(1337,5604)
	(1364,5624)(1393,5645)(1423,5667)
	(1455,5690)(1488,5713)(1522,5738)
	(1557,5762)(1593,5787)(1628,5812)
	(1662,5835)(1694,5857)(1723,5877)
	(1748,5895)(1770,5909)(1786,5920)
	(1798,5929)(1806,5934)(1810,5937)(1812,5938)
\path(4212,2788)(4213,2789)(4216,2791)
	(4221,2794)(4228,2800)(4239,2808)
	(4254,2818)(4274,2832)(4297,2849)
	(4326,2869)(4359,2893)(4396,2920)
	(4438,2950)(4484,2983)(4534,3019)
	(4588,3057)(4644,3097)(4703,3139)
	(4764,3182)(4826,3227)(4890,3272)
	(4954,3317)(5018,3363)(5082,3408)
	(5145,3453)(5207,3497)(5269,3540)
	(5329,3582)(5387,3623)(5444,3663)
	(5499,3702)(5553,3739)(5605,3776)
	(5656,3810)(5704,3844)(5752,3877)
	(5797,3908)(5842,3938)(5885,3967)
	(5926,3995)(5967,4022)(6007,4049)
	(6045,4074)(6083,4099)(6120,4124)
	(6157,4147)(6193,4171)(6229,4193)
	(6264,4216)(6299,4238)(6300,4238)
	(6341,4264)(6383,4290)(6424,4315)
	(6465,4340)(6507,4365)(6549,4390)
	(6590,4415)(6633,4440)(6675,4464)
	(6717,4489)(6760,4513)(6803,4537)
	(6846,4561)(6890,4585)(6933,4609)
	(6977,4633)(7020,4656)(7064,4679)
	(7108,4702)(7151,4724)(7195,4746)
	(7238,4768)(7281,4789)(7323,4810)
	(7365,4830)(7407,4850)(7448,4870)
	(7489,4888)(7529,4907)(7569,4925)
	(7607,4942)(7646,4959)(7683,4975)
	(7721,4991)(7757,5006)(7793,5021)
	(7828,5035)(7863,5049)(7898,5062)
	(7932,5075)(7966,5088)(8000,5100)
	(8037,5114)(8074,5127)(8111,5140)
	(8148,5152)(8185,5165)(8224,5177)
	(8263,5189)(8302,5201)(8343,5213)
	(8385,5225)(8429,5238)(8474,5250)
	(8520,5263)(8569,5275)(8619,5288)
	(8670,5302)(8724,5315)(8779,5328)
	(8835,5342)(8891,5356)(8948,5370)
	(9005,5383)(9062,5396)(9117,5409)
	(9169,5421)(9219,5433)(9265,5443)
	(9307,5453)(9343,5461)(9375,5469)
	(9402,5474)(9422,5479)(9438,5483)
	(9449,5485)(9457,5487)(9460,5488)(9462,5488)
\path(9462,5488)(10062,5563)
\path(10062,5563)(10064,5562)(10069,5561)
	(10078,5558)(10091,5554)(10110,5548)
	(10134,5541)(10163,5532)(10196,5521)
	(10232,5510)(10271,5498)(10311,5485)
	(10352,5472)(10393,5458)(10433,5445)
	(10472,5433)(10509,5420)(10545,5408)
	(10579,5397)(10611,5386)(10642,5375)
	(10672,5364)(10701,5354)(10729,5343)
	(10756,5333)(10783,5322)(10810,5311)
	(10837,5301)(10837,5300)(10862,5290)
	(10888,5279)(10914,5268)(10941,5257)
	(10968,5245)(10996,5232)(11026,5218)
	(11057,5204)(11090,5189)(11125,5173)
	(11161,5156)(11199,5138)(11239,5118)
	(11281,5098)(11324,5078)(11368,5056)
	(11413,5035)(11457,5013)(11500,4992)
	(11540,4972)(11578,4954)(11612,4937)
	(11641,4923)(11665,4911)(11683,4902)
	(11697,4896)(11705,4891)(11710,4889)(11712,4888)
\put(3987,2863){\makebox(0,0)[lb]{\smash{{{\SetFigFont{6}{7.2}{rm}$\omega$}}}}}
\put(2787,6463){\makebox(0,0)[lb]{\smash{{{\SetFigFont{6}{7.2}{rm}$g_{j+2}(x)$}}}}}
\put(10887,5938){\makebox(0,0)[lb]{\smash{{{\SetFigFont{6}{7.2}{rm}$g_j(x)$}}}}}
\put(3837,4363){\makebox(0,0)[lb]{\smash{{{\SetFigFont{6}{7.2}{rm}$g_{j+1}(x)$}}}}}
\put(9687,1813){\makebox(0,0)[lb]{\smash{{{\SetFigFont{6}{7.2}{rm}$g_{j-1}(x)$}}}}}
\put(1587,6163){\makebox(0,0)[lb]{\smash{{{\SetFigFont{6}{7.2}{rm}$z_{i-1}$}}}}}
\put(5937,2038){\makebox(0,0)[lb]{\smash{{{\SetFigFont{6}{7.2}{rm}$R_i$}}}}}
\put(9387,5638){\makebox(0,0)[lb]{\smash{{{\SetFigFont{6}{7.2}{rm}$z_i$}}}}}
\put(2487,5338){\makebox(0,0)[lb]{\smash{{{\SetFigFont{6}{7.2}{rm}$\psi_i(x)$}}}}}
\put(1737,2443){\makebox(0,0)[lb]{\smash{{{\SetFigFont{6}{7.2}{rm}$x_{i-1}$}}}}}
\put(4437,2458){\makebox(0,0)[lb]{\smash{{{\SetFigFont{6}{7.2}{rm}$y_i^-$}}}}}
\put(5487,2458){\makebox(0,0)[lb]{\smash{{{\SetFigFont{6}{7.2}{rm}$y_i^+$}}}}}
\put(9387,2443){\makebox(0,0)[lb]{\smash{{{\SetFigFont{6}{7.2}{rm}$x_i$}}}}}
\end{picture}
}